\documentclass[leqno,draft]{article}

\def\re{\mathop{\rm Re}}
\def\im{\mathop{\rm Im}}

\newtheorem{theorem}{Theorem}
\newtheorem{lemma}[theorem]{Lemma}
\newtheorem{proposition}[theorem]{Proposition}
\newtheorem{definition}[theorem]{Definition}
\newtheorem{corollary}[theorem]{Corollary}

\newcommand{\begintheorem}{\addtocounter{equation}{1}\begin{theorem}}
\newcommand{\beginlemma}{\addtocounter{equation}{1}\begin{lemma}}
\newcommand{\beginproposition}{\addtocounter{equation}{1}\begin{proposition}}
\newcommand{\begindefinition}{\addtocounter{equation}{1}\begin{definition}}
\newcommand{\begincorollary}{\addtocounter{equation}{1}\begin{corollary}}

\begin{document}

\title{Sums and averages}

\author{Stephen Semmes \\
        Rice University}

\date{}

\maketitle

\begin{abstract}
These informal notes are concerned with sums and averages in various
situations in analysis.
\end{abstract}

\tableofcontents

\section{Real and complex numbers}
\label{real and complex numbers}
\setcounter{equation}{0}

        Let ${\bf R}$ be the field of real numbers.  The \emph{absolute value}
of a real number $x$ is denoted $|x|$ and defined to be equal to $x$ when
$x \ge 0$ and to $-x$ when $x \le 0$.  Thus $|x| \ge 0$,
\begin{equation}
        |x + y| \le |x| + |y|,
\end{equation}
and
\begin{equation}
        |x \, y| = |x| \, |y|
\end{equation}
for every $x, y \in {\bf R}$.

        If $A$ is a nonempty set of real numbers with an upper bound,
then there is a unique real number which is the least upper bound or
supremum of $A$, denoted $\sup A$.  Similarly, a nonempty set $A$ of
real numbers with a lower bound has a greatest lower bound or infimum,
denoted $\inf A$.  It is sometimes convenient to put $\sup A = +
\infty$ when $A$ has no upper bound, or $\inf A = -\infty$ when $A$
has no lower bound.  Normally we shall only be concerned with infima
of sets of nonnegative real numbers here, which are also nonnegative.

        Let ${\bf C}$ be the field of complex numbers.  A complex
number $z$ can be expressed as $x + y \, i$, where $x, y \in {\bf R}$
and $i^2 = -1$.  In this case, $x$ and $y$ are known as the real and
imaginary parts of $z$, and are denoted $\re z$, $\im z$,
respectively.  The \emph{complex conjugate} $\overline{z}$ of $z$ is
defined by
\begin{equation}
        \overline{z} = x - y \, i.
\end{equation}
In particular, $\overline{\overline{z}} = z$ and
\begin{equation}
        \re z = \frac{z + \overline{z}}{2},
         \quad \im z = \frac{z - \overline{z}}{2 \, i}.
\end{equation}
Observe that
\begin{equation}
        \overline{z + w} = \overline{z} + \overline{w}
\end{equation}
and
\begin{equation}
        \overline{z \, w} = \overline{z} \, \overline{w}
\end{equation}
for every $z, w \in {\bf C}$.  The \emph{modulus} of $z = x + y \, i$
is the nonnegative real number given by
\begin{equation}
        |z| = (x^2 + y^2)^{1/2}.
\end{equation}
Equivalently,
\begin{equation}
        |z|^2 = z \, \overline{z},
\end{equation}
and hence
\begin{equation}
        |z \, w| = |z| \, |w|
\end{equation}
for every $z, w \in {\bf C}$.

        Of course, $|\re z|, |\im z| \le |z|$ for every $z \in {\bf
C}$.  If $z, w \in {\bf C}$, then
\begin{eqnarray}
         |z + w|^2 & = & (z + w) (\overline{z} + \overline{w}) \\
        & = & |z|^2 + z \, \overline{w} + \overline{z} \, w + |w|^2 \nonumber\\
        & = & |z|^2 + 2 \, \re z \, \overline{w} + |w|^2, \nonumber \\
\end{eqnarray}
since $\overline{z \, \overline{w}} = \overline{z} \, w$.  This implies
that
\begin{equation}
        |z + w|^2 \le |z|^2 + 2 \, |z| \, |w| + |w|^2 = (|z| + |w|)^2,
\end{equation}
and therefore
\begin{equation}
        |z + w| \le |z| + |w|.
\end{equation}

\section{Cesaro means}
\label{cesaro means}
\setcounter{equation}{0}

        As usual, a sequence $\{z_j\}_{j = 0}^\infty$ of complex
numbers converges to $z \in {\bf C}$ if for every $\epsilon > 0$ there
is a nonnegative integer $L$ such that
\begin{equation}
        |z_j - z| < \epsilon
\end{equation}
for every $j \ge L$.  In this case, one can show that the sequence of averages
\begin{equation}
        \zeta_n = \frac{z_1 + \cdots + z_n}{n + 1}
\end{equation}
also converges to $z$ as $n \to \infty$.  However, there are also
sequences $\{z_j\}_{j = 0}^\infty$ of complex numbers that do not
converge, but for which the corresponding sequence $\{\zeta_n\}_{n =
0}^\infty$ of averages does converge.  For example, if $z_j = (-1)^j$,
then $\zeta_n = 0$ when $n$ is odd and $\zeta_n = 1/(n + 1)$ when $n$
is even, and $\lim_{n \to \infty} \zeta_n = 0$.

        If $z$ is any complex number, then
\begin{equation}
        (z - 1) \sum_{j = 0}^n z^j = z^{n + 1} - 1,
\end{equation}
where $z^j = 1$ for every $z \in {\bf C}$ when $j = 0$.  Hence
\begin{equation}
        \sum_{j = 0}^n z^j = \frac{z^{n + 1} - 1}{z - 1}
\end{equation}
when $z \ne 1$.  It follows that
\begin{equation}
        \lim_{n \to \infty} \frac{1}{n + 1} \sum_{j = 0}^n z^j = 0
\end{equation}
when $|z| = 1$ and $z \ne 1$.  This extends the case of $z = -1$
described in the previous paragraph.

        Let $\sum_{j = 0}^\infty a_j$ be an infinite series of complex
numbers, and consider the sequence of partial sums
\begin{equation}
        b_l = \sum_{j = 0}^l a_j.
\end{equation}
By definition, $\sum_{j = 0}^\infty a_j$ converges if $\{b_l\}_{l =
0}^\infty$ converges, in which event
\begin{equation}
        \sum_{j = 0}^\infty a_j = \lim_{l \to \infty} b_l.
\end{equation}
The average $\beta_n$ of $b_0, \ldots, b_n$ is also given by
\begin{equation}
        \beta_n = \sum_{j = 0}^n \frac{n + 1 - j}{n + 1} \, a_j.
\end{equation}
The series $\sum_{j = 0}^\infty a_j$ is said to be \emph{Cesaro
summable} if $\{\beta_n\}_{n = 0}^\infty$ converges.

        Consider the case of a geometric series
\begin{equation}
        \sum_{j = 0}^\infty a^j,
\end{equation}
where $a \in {\bf C}$ and $a^j = 1$ when $j = 0$ again.  If $|a| < 1$,
then $\lim_{j \to \infty} a^j = 0$, and the corresponding geometric
series converges with
\begin{equation}
        \sum_{j = 0}^\infty a^j = \frac{1}{1 - a}
\end{equation}
by the previous computations.  If $|a| \ge 1$, then $|a^j| = |a|^j \ge
1$ for every $j$, and the geometric series does not converge in the
conventional sense.  It is Cesaro summable with sum $1 / (1 - a)$ when
$|a| = 1$ and $a \ne 1$, by the earlier computations for the partial
sums applied twice to estimate their averages too.

\section{Admissible series}
\label{admissible series}
\setcounter{equation}{0}

        Let $\sum_{j = 0}^\infty a_j$ be an infinite series of complex
numbers.  If
\begin{equation}
        \sum_{j = 0}^\infty a_j \, t^j
\end{equation}
converges for some real number $t \le 1$, then
\begin{equation}
        \lim_{j \to \infty} a_j \, t^j = 0,
\end{equation}
which implies that $\{a_j \, t^j\}_{j = 0}^\infty$ is a bounded sequence.
Conversely, if $\{a_j \, t^j\}_{j = 0}^\infty$ is bounded, then
\begin{equation}
        \sum_{j = 0}^\infty |a_j| \, r^j
\end{equation}
converges for $0 \le r < t$, by comparison with a convergent geometric series.

        Let us say that an infinite series $\sum_{j = 0}^\infty a_j$
of complex numbers is \emph{admissible} if any of the previous
conditions holds for every positive real number $r < 1$ or $t < 1$, as
appropriate, so that each of the other conditions also holds for every
$0 \le r < 1$ or $0 \le t < 1$.  This is the same as saying that
\begin{equation}
        f(z) = \sum_{j = 0}^\infty a_j \, z^j
\end{equation}
has radius of convergence greater than or equal to $1$.  If $b_l =
\sum_{j = 0}^l a_j$ are the partial sums of $\sum_{j = 0}^\infty a_j$,
then $\sum_{l = 0}^\infty b_l$ is admissible too.  Indeed, $a_j =
O(R^j)$ implies that $b_l = O(R^l)$ for each $R > 1$.

        Put $b_{-1} = 0$, so that
\begin{equation}
        \sum_{j = 0}^n a_j \, z^j = \sum_{j = 0}^n (b_j - b_{j - 1}) \, z^j
 = \sum_{j = 0}^n b_j \, z^j - \sum_{j = 0}^n b_{j - 1} \, z^j.
\end{equation}
Since
\begin{equation}
 \sum_{j = 0}^n b_{j - 1} \, z^j = \sum_{j = 1}^n b_{j - 1} \, z^j
                                  = \sum_{j = 0}^{n - 1} b_j \, z^{j + 1},
\end{equation}
we get that
\begin{equation}
        \sum_{j = 0}^n a_j \, z^j
 = \sum_{j = 0}^{n - 1} b_j \, (z^j - z^{j + 1}) + b_n \, z^n
 = (1 - z) \sum_{j = 0}^{n - 1} b_j \, z^j + b_n \, z^n.
\end{equation}
If $|z| < 1$, then admissibility of $\sum_{n = 0}^\infty b_n$ implies
that $\lim_{n \to \infty} b_n \, z^n = 0$, and
\begin{equation}
        f(z) = (1 - z) \sum_{j = 0}^\infty b_j \, z^j.
\end{equation}

\section{Abel summability}
\label{abel summability}
\setcounter{equation}{0}

        If $\sum_{j = 0}^\infty a_j$ converges, then it is well known that
\begin{equation}
        \lim_{r \to 1-} \sum_{j = 0}^\infty a_j \, r^j
\end{equation}
exists and is equal to $\sum_{j = 0}^\infty a_j$.  An admissible
series $\sum_{j = 0}^\infty a_j$ of complex numbers is said to be
\emph{Abel summable} when this limit exists.

        For example, $\sum_{j = 0}^\infty a^j$ is admissible for any
complex number $a$ with $|a| = 1$, and
\begin{equation}
        \sum_{j = 0}^\infty a^j \, z^j = \frac{1}{1 - a \, z}
\end{equation}
for every $z \in {\bf C}$ with $|z| < 1$, which implies that $\sum_{j
= 0}^\infty a^j$ is Abel summable to $1/(1 - a)$ when $a \ne 1$.
Similarly, $\sum_{j = 0}^\infty (j + 1) \, a^j$ is admissible when
$|a| = 1$, and
\begin{equation}
        \sum_{j = 0}^\infty (j + 1) \, a^j \, z^j = \frac{1}{(1 - a \, z)^2}
\end{equation}
for $|z| < 1$, so that $\sum_{j = 0}^\infty j \, a^j$ is Abel summable
to $1/(1 - a)^2$ when $a \ne 1$.

        Let $\{b_j\}_{j = 0}^\infty$ be a sequence of complex numbers,
and consider
\begin{equation}
        \beta_n = \frac{1}{n + 1} \sum_{j = 0}^n b_j.
\end{equation}
If $\{\beta_n\}_{n = 0}^\infty$ converges, then
\begin{equation}
 \frac{n}{n + 1} \, \beta_{n - 1} = \frac{1}{n + 1} \sum_{j = 0}^{n - 1} b_j
\end{equation}
converges to the same value.  This implies that
\begin{equation}
        \frac{b_n}{n + 1} = \beta_n - \frac{n}{n + 1} \, \beta_{n - 1} \to 0
\end{equation}
as $n \to \infty$.  If $\sum_{j = 0}^\infty a_j$ is a Cesaro summable
series of complex numbers, then one cam apply this to $b_n = \sum_{j =
0}^n a_j$ to get that
\begin{equation}
        \frac{a_n}{n + 1} = \frac{b_n - b_{n - 1}}{n + 1} \to 0
\end{equation}
as $n \to \infty$.  In particular, $\sum_{j = 0}^\infty (j + 1) \,
a^j$ is not Cesaro summable for any $a \in {\bf C}$ with $|a| = 1$.

        If $\sum_{j = 0}^\infty a_j$ is Cesaro summable, then $\sum_{j
= 0}^\infty a_j$ is an admissible series, since $a_n = O(n)$.  It is
well known that $\sum_{j = 0}^\infty a_j$ is Abel summable in this
case, and with the same sum.  For if $b_l = \sum_{j = 0}^l a_j$ are
the partial sums and $c_n = \sum_{l = 0}^n b_l$ are their partial sums, then
\begin{equation}
 \sum_{j = 0}^\infty a_j \, z^j = (1 - z) \sum_{j = 0}^\infty b_j \, z^j
 = (1 - z)^2 \sum_{j = 0}^\infty c_j \, z^j
\end{equation}
when $|z| < 1$, as in the preceding section.  Equivalently,
\begin{equation}
        \sum_{j = 0}^\infty a_j \, z^j
 = (1 - z)^2 \sum_{j = 0}^\infty (j + 1) \, \beta_j \, z^j,
\end{equation}
when $|z| < 1$, where $\beta_j = c_j/(j + 1)$ is the average of $b_0,
\ldots, b_j$, as before.  Of course,
\begin{equation}
        \sum_{j = 0}^\infty (j + 1) \, z^j = \frac{1}{(1 - z)^2}
\end{equation}
is the derivative of the usual geometric series.  If the $\beta_j$'s
were constant, then the desired conclusion would follow immediately.
If $\beta_j \to 0$ as $j \to \infty$, then one can check that $\sum_{j
= 0}^\infty a_j$ is Abel summable with sum equal to $0$.  If $\sum_{j
= 0}^\infty a_j$ is Cesaro summable, so that $\{\beta_j\}_{j =
0}^\infty$ converges, then one can combine these two cases to show
that $\sum_{j = 0}^\infty a_j$ is Abel summable with the same sum.

\section{Cauchy products}
\label{cauchy products}
\setcounter{equation}{0}

        If $\sum_{j = 0}^\infty a_j$ and $\sum_{l = 0}^\infty b_l$ are
infinite series of complex numbers, then their \emph{Cauchy product}
is the infinite series $\sum_{n = 0}^\infty c_n$ whose terms are given
by
\begin{equation}
        c_n = \sum_{j = 0}^n a_j \, b_{n - j}.
\end{equation}
Formally,
\begin{equation}
\label{sum c_n = (sum a_j) (sum b_l)}
        \sum_{n = 0}^\infty c_n =
 \Big(\sum_{j = 0}^\infty a_j\Big) \Big(\sum_{l = 0}^\infty b_l\Big),
\end{equation}
and in particular this holds when $a_j = b_l = 0$ for all but finitely
many $j$ and $l$.  Moreover, (\ref{sum c_n = (sum a_j) (sum b_l)})
holds when $a_j = 0$ for all but finitely many $j$ or $b_l = 0$ for
all but finitely many $l$.  Note that $\sum_{n = 0}^\infty c_n \, z^n$
is the Cauchy product of $\sum_{j = 0}^\infty a_j \, z^j$ and $\sum_{l
= 0}^\infty b_l \, z^l$ for every $z \in {\bf C}$.

        If $\sum_{j = 0}^\infty a_j$ and $\sum_{l = 0}^\infty b_l$
converge absolutely, which means that $\sum_{j = 0}^\infty |a_j|$ and
$\sum_{l = 0}^\infty |b_l|$ converge, then $\sum_{n = 0}^\infty c_n$
converges absolutely, and (\ref{sum c_n = (sum a_j) (sum b_l)}) holds.
In connection with this, observe that
\begin{equation}
        |c_n| \le \sum_{j = 0}^n |a_j| \, |b_{n - j}|,
\end{equation}
where the sum on the right side of the inequality corresponds exactly
to the Cauchy product of $\sum_{j = 0}^\infty |a_j|$ and $\sum_{l =
0}^\infty |b_l|$.  One can also show that $\sum_{n = 0}^\infty c_n$
converges and satisfies (\ref{sum c_n = (sum a_j) (sum b_l)}) when
$\sum_{j = 0}^\infty a_j$ and $\sum_{l = 0}^\infty b_l$ converge and
at least one of these two series converges absolutely.

        If $a_j, b_j = O(R^j)$ for some $R > 0$, then $c_n = O(n \,
R^n)$.  This implies that $\sum_{n = 0}^\infty c_n$ is admissible when
$\sum_{j = 0}^\infty a_j$ and $\sum_{l = 0}^\infty b_l$ are
admissible.  In this case, the corresponding power series $\sum_{j =
0}^\infty a_j \, z^j$, $\sum_{l = 0}^\infty b_l \, z^l$, and $\sum_{n
= 0}^\infty c_n \, z^n$ converge absolutely for every $z \in {\bf C}$
with $|z| < 1$, and satisfy
\begin{equation}
\label{sum c_n z^n = (sum a_j z^j) (sum b_l z^l)}
 \sum_{n = 0}^\infty c_n \, z^n = \Big(\sum_{j = 0}^\infty a_j \, z^j\Big)
                                   \Big(\sum_{l = 0}^\infty b_l \, z^l \Big).
\end{equation}

        If $\sum_{j = 0}^\infty a_j$ and $\sum_{l = 0}^\infty b_l$ are
Abel summable, then $\sum_{n = 0}^\infty c_n$ is Abel summable, and
their Abel sums satisfy (\ref{sum c_n = (sum a_j) (sum b_l)}).  This
follows from the analogous statement for the corresponding power
series on the unit disk, as in the previous paragraph.

\section{Norms on vector spaces}
\label{norms on vector spaces}
\setcounter{equation}{0}

        Let $V$ be a real or complex vector space.  A real-valued
function $N$ on $V$ is said to be a \emph{norm} if $N(v) \ge 0$
for every $v \in V$, $N(v) = 0$ if and only if $v = 0$,
\begin{equation}
        N(t \, v) = |t| \, N(v)
\end{equation}
for every $v \in V$ and real or complex number $t$, as appropriate, and
\begin{equation}
        N(v + w) \le N(v) + N(w)
\end{equation}
for every $v, w \in V$.  Thus the absolute value and modulus determine
norms on ${\bf R}$ and ${\bf C}$ as one-dimensional vector spaces,
respectively.

        Let $n$ be a positive integer, and consider the spaces ${\bf
R}^n$ and ${\bf C}^n$ of $n$-tuples of real and complex numbers.  As
usual, these are vector spaces with respect to coordinatewise addition
and multiplication.  Consider
\begin{equation}
        \|v\|_1 = |v_1| + \cdots + |v_n|
\end{equation}
and
\begin{equation}
        \|v\|_\infty = \max(|v_1|, \ldots, |v_n|)
\end{equation}
for $v = (v_1, \ldots, v_n)$ in ${\bf R}^n$ or ${\bf C}^n$.  It is
easy to see that $\|v\|_1$ and $\|v\|_\infty$ are norms on ${\bf R}^n$
and ${\bf C}^n$.

        The standard Euclidean norm on ${\bf R}^n$ and ${\bf C}^n$ is
defined by
\begin{equation}
        \|v\|_2 = (|v_1|^2 + \cdots + |v_n|^2)^{1/2}.
\end{equation}
This clearly satisfies the positivity and homogeneity requirements of
a norm.  The triangle inequality will be discussed in the next two sections.

        Observe that
\begin{equation}
        \|v\|_\infty \le \|v\|_1
\end{equation}
and
\begin{equation}
        \|v\|_\infty \le \|v\|_2
\end{equation}
for every $v$ in ${\bf R}^n$ or ${\bf C}^n$.  One can also check that
\begin{equation}
        \|v\|_2 \le \|v\|_1
\end{equation}
for every $v$ in ${\bf R}^n$ or ${\bf C}^n$, using the first inequality.
More precisely,
\begin{equation}
        \|v\|_2^2 \le \|v\|_1 \, \|v\|_\infty \le \|v\|_1^2,
\end{equation}
where the first step follows directly from the definitions.

\section{Inner product spaces}
\label{inner product spaces}
\setcounter{equation}{0}

        Let $V$ be a real or complex vector space again.  An
\emph{inner product} on $V$ is a real or complex-valued function
$\langle v, w \rangle$ defined for $v, w \in V$ with the following
properties.  First, $\langle v, w \rangle$ is linear as a function of
$v$ for each fixed $w \in V$.  Second,
\begin{equation}
        \langle w, v \rangle = \langle v, w \rangle
\end{equation}
for every $v, w \in V$ in the real case, and
\begin{equation}
        \langle w, v \rangle = \overline{\langle v, w \rangle}
\end{equation}
for every $v, w \in V$ in the complex case.  This implies that
$\langle v, w \rangle$ is linear in $w$ in the real case, and that it
is conjugate-linear in the complex case.  This also implies that
$\langle v, v \rangle \in {\bf R}$ for every $v \in V$ in the complex
case.  The third condition asks that
\begin{equation}
        \langle v, v \rangle > 0
\end{equation}
for every $v \in V$ with $v \ne 0$.  Of course, $\langle 0, 0 \rangle
= 0$ by linearity.

        For each $v \in V$, let $\|v\|$ be the nonnegative real number
defined by
\begin{equation}
        \|v\| = \langle v, v \rangle^{1/2}.
\end{equation}
The \emph{Cauchy--Schwarz inequality} states that
\begin{equation}
        |\langle v, w \rangle| \le \|v\| \, \|w\|
\end{equation}
for every $v, w \in V$.  To show this, one can start with
\begin{equation}
        \langle v + t \, w, v + t \, w \rangle \ge 0
\end{equation}
for every $t \in {\bf R}$ or ${\bf C}$, as appropriate.  This implies that
\begin{equation}
        2 \, |t| \, |\langle v, w \rangle| \le \|v\|^2 + |t|^2 \, \|w\|^2,
\end{equation}
by expanding the inner product and collecting terms, and choosing the
sign of $t$ in the real case or $t/|t|$ in the complex case to get the
absolute value or modulus of the inner product on the left.  The
Cauchy--Schwarz inequality follows by taking $|t| = \|v\| / \|w\|$
when $v, w \ne 0$.

        Using the Cauchy--Schwarz inequality, one gets that
\begin{eqnarray}
        \|v + w\|^2 & = & \langle v + w, v + w \rangle \\
 & \le & \|v\|^2 + 2 \, \|v\| \, \|w\| + \|w\|^2 = (\|v\| + \|w\|)^2.\nonumber
\end{eqnarray}
Hence
\begin{equation}
        \|v + w\| \le \|v\| + \|w\|
\end{equation}
for every $v, w \in V$.  Thus $\|v\|$ is a norm on $V$, since the
positivity and homogeneity conditions are clearly satisfied.

        The standard inner products on ${\bf R}^n$, ${\bf C}^n$ are given by
\begin{equation}
        \langle v, w \rangle = \sum_{j = 1}^n v_j \, w_j
\end{equation}
for $v = (v_1, \ldots, v_n), w = (w_1, \ldots, w_n) \in {\bf R}^n$, and
\begin{equation}
        \langle v, w \rangle = \sum_{j = 1}^n v_j \, \overline{w_j}
\end{equation}
in the complex case.  The corresponding norm
\begin{equation}
        \|v\| = \Big(\sum_{j = 1}^n |v_j|^2\Big)^{1/2}
\end{equation}
is the same as the standard Euclidean norm $\|v\|_2$.

\section{Convexity}
\label{convexity}
\setcounter{equation}{0}

        A set $E$ in a real or complex vector space $V$ is said to be
\emph{convex} if
\begin{equation}
\label{t v + (1 - t) w in E}
        t \, v + (1 - t) \, w \in E
\end{equation}
for every $v, w \in E$ and real number $t$ such that $0 < t < 1$.
For example, if $N$ is a norm on $V$, then the closed unit ball
\begin{equation}
        B = \{v \in V : N(v) \le 1\}
\end{equation}
is convex.

        Conversely, if $N$ is a real-valued function on $V$ which
satisfies the positivity and homogeneity conditions of a norm, and if
the corresponding closed unit ball $B$ is convex, then one can show
that $N$ satisfies the triangle inequality and hence is a norm.  For
if $\widehat{v}$, $\widehat{w}$ are nonzero vectors in $V$, then we
can apply (\ref{t v + (1 - t) w in E}) with $v = \widehat{v} /
N(\widehat{v})$, $w = \widehat{w} / N(\widehat{w})$,
\begin{equation}
        t = \frac{N(\widehat{v})}{N(\widehat{v}) + N(\widehat{w})},
\end{equation}
and $E = B$ to get that
\begin{equation}
 \frac{\widehat{v} + \widehat{w}}{N(\widehat{v}) + N(\widehat{w})} \in B,
\end{equation}
which says exactly that
\begin{equation}
        N(\widehat{v} + \widehat{w}) \le N(\widehat{v}) + N(\widehat{w}),
\end{equation}
as desired.

        A real-valued function $\phi$ on the real line is convex if
\begin{equation}
        \phi(t \, x + (1 - t) \, y) \le t \, \phi(x) + (1 - t) \, \phi(y)
\end{equation}
for every $x, y, t \in {\bf R}$ with $0 < t < 1$.  For example,
\begin{equation}
        \phi(x) = |x|^p
\end{equation}
is convex when $p \ge 1$.

        If $p \ge 1$ and $v \in {\bf R}^n$ or ${\bf C}^n$, then put
\begin{equation}
        \|v\|_p = \Big(\sum_{j = 1}^n |v_j|^p\Big)^{1/p}.
\end{equation}
This satisfies the positivity and homogeneity requirements of a norm,
and one can use the convexity of $\phi(x) = |x|^p$ to show that the
corresponding unit ball is convex, and hence that $\|v\|_p$ is a norm.

\section{A few estimates}
\label{a few estimates}
\setcounter{equation}{0}

        Observe that
\begin{equation}
        \|v\|_\infty \le \|v\|_p
\end{equation}
for every $v \in {\bf R}^n$, ${\bf C}$ and $p \ge 1$.  If $1 \le p < q
< \infty$, then
\begin{equation}
        \|v\|_q^q \le \|v\|_p^p \, \|v\|_\infty^{q - p} \le \|v\|_p^q,
\end{equation}
which implies that
\begin{equation}
        \|v\|_q \le \|v\|_p.
\end{equation}
In the other direction, it is easy to see that
\begin{equation}
        \|v\|_p \le n^{1/p} \, \|v\|_\infty
\end{equation}
for every $v \in {\bf R}^n$ or ${\bf C}^n$ and $p \ge 1$.

        If $1 \le p < q < \infty$, then
\begin{equation}
        \|v\|_p \le n^{1/p - 1/q} \, \|v\|_q
\end{equation}
for every $v \in {\bf R}^n$ or ${\bf C}^n$.  Equivalently,
\begin{equation}
        \Big(\frac{1}{n}\sum_{j = 1}^n |v_j|^p\Big)^{1/p}
         \le \Big(\frac{1}{n}\sum_{j = 1}^n |v_j|^q\Big)^{1/q}.
\end{equation}
More precisely, this is the same as
\begin{equation}
        \Big(\frac{1}{n}\sum_{j = 1}^n |v_j|^p\Big)^{q/p}
         \le \frac{1}{n}\sum_{j = 1}^n |v_j|^q,
\end{equation}
which can be derived from the convexity of $\phi(x) = |x|^{q/p}$ on
the real line.

        If $N$ is any norm on ${\bf R}^n$ or ${\bf C}^n$, then
\begin{equation}
        N(v) \le A \, \|v\|_2
\end{equation}
for some $A > 0$ and every $v \in {\bf R}^n$ or ${\bf C}^n$, as appropriate.
This can be verified using the triangle inequality to estimate $N(v)$ in terms
of the norms of the standard basis vectors.  The triangle inequality also
implies that
\begin{equation}
        N(v) - N(w), N(w) - N(v) \le N(v - w)
\end{equation}
for every $v, w \in {\bf R}^n$ or ${\bf C}^n$, and hence that
\begin{equation}
        |N(v) - N(w)| \le N(v - w) \le A \, \|v - w\|_2.
\end{equation}
This shows that $N$ is continuous with respect to the standard
topology on ${\bf R}^n$ or ${\bf C}^n$.

        The unit sphere in ${\bf R}^n$ or ${\bf C}^n$ with respect to
the standard Euclidean norm consists of the vectors $v$ such that
$\|v\|_2 = 1$.  By compactness, the continuous function $N$ attains
its minimum on the unit sphere, which is positive.  Thus
\begin{equation}
        N(v) \ge a
\end{equation}
for some $a > 0$ and every $v \in {\bf R}^n$ or ${\bf C}^n$ with $\|v\|_2 = 1$.
It follows that
\begin{equation}
        a \, \|v\|_2 \le N(v)
\end{equation}
for every $v \in {\bf R}^n$ or ${\bf C}^n$, by homogeneity.

\section{Operator norms}
\label{operator norms}
\setcounter{equation}{0}

        Let $V$, $W$ be vector spaces, both real or both complex, with
norms $\|\cdot \|_V$, $\|\cdot \|_W$, respectively.  A linear mapping
$T : V \to W$ is said to be \emph{bounded} if there is an $A \ge 0$ such that
\begin{equation}
\label{||T(v)||_W le A ||v||_V}
        \|T(v)\|_W \le A \, \|v\|_V
\end{equation}
for every $v \in V$.  If $V$ is ${\bf R}^n$ or ${\bf C}^n$ and
$\|v\|_V = \|v\|_p$ for some $p$, $1 \le p \le \infty$, then every
linear mapping $T : V \to W$ is bounded, as one can see by expressing
any $v \in V$ as a linear combination of the standard basis vectors.
This also works for any norm on ${\bf R}^n$ or ${\bf C}^n$, since any
norm is equivalent to the $p$-norms, as in the previous section.  The
same statement holds as well for any norm on any finite-dimensional
vector space $V$, because $V$ is then isomorphic to ${\bf R}^n$ or
${\bf C}^n$ for some $n$.

        The \emph{operator norm} of a bounded linear mapping $T : V
\to W$ is defined by
\begin{equation}
        \|T\|_{op} = \sup \{\|T(v)\|_W : v \in V, \, \|v\|_V \le 1\}.
\end{equation}
Equivalently, (\ref{||T(v)||_W le A ||v||_V}) holds with $A =
\|T\|_{op}$, and $\|T\|_{op}$ is the smallest nonnegative real number
with this property.  One can check that the operator norm is a norm on
the vector space of bounded linear mappings from $V$ into $W$.

        Suppose that $V_1$, $V_2$, and $V_3$ are vector spaces, all
real or all complex, equipped with norms as before, and that $T_1 :
V_1 \to V_2$ and $T_2 : V_2 \to V_3$ are bounded linear mappings.
It is easy to see that the composition $T_2 \circ T_1$, defined by
\begin{equation}
        (T_2 \circ T_1)(v) = T_2(T_1(v)), \ v \in V_1,
\end{equation}
is a bounded linear mapping from $V_1$ into $V_3$, and that
\begin{equation}
        \|T_2 \circ T_1\|_{op, 13} \le \|T_1\|_{op, 12} \, \|T_2\|_{op, 23},
\end{equation}
where the subscripts indicate the spaces involved in the operator norms.
For any normed vector space $V$, the identity operator $I$ defined by
$I(v) = v$ for each $v \in V$ is a bounded linear mapping from $V$ into
itself, and satisfies
\begin{equation}
        \|I\|_{op} = 1,
\end{equation}
using the same norm on $V$ as both the domain and range.

        Suppose that $V$ is ${\bf R}^n$ or ${\bf C}^n$ equipped with
the norm
\begin{equation}
        \|v\|_1 = \sum_{j = 1}^n |v_j|,
\end{equation}
and let $e_1, \ldots, e_n$ be the standard basis of $V$, so that the $j$th
coordinate of $e_j$ is equal to $1$ and the other coordinates are $0$.
In this case,
\begin{equation}
        \|T\|_{op} = \max \{\|T(e_1)\|_W, \ldots, \|T(e_n)\|_W\}
\end{equation}
for any linear mapping $T : V \to W$.  If $V$ is any vector space with
any norm $\|\cdot \|_V$ and the norm $\|\cdot \|_W$ on $W$ is
associated to an inner product $\langle \cdot, \cdot \rangle_W$, then
\begin{equation}
        \|T\|_{op} = \sup \{|\langle T(v), w\rangle_W| :
                     v \in V, \, \|v\|_V \le 1, \, w \in W, \, \|w\|_W \le 1\}.
\end{equation}
Indeed, for each $z \in W$,
\begin{equation}
  \|z\|_W = \sup \{|\langle z, w \rangle_W| : w \in W, \, \|w\|_W \le 1\},
\end{equation}
since the inner product is bounded by the norm because of the
Cauchy--Schwarz inequality, and equality occurs with $w = z / \|z\|_W$
when $z \ne 0$ and with any $w$ when $z = 0$.

\section{Linear mappings}
\label{linear mappings}
\setcounter{equation}{0}

        Let $T$ be a linear mapping from a vector space $V$ into
itself, and let $T^j$ be the composition of $j$ factors of $T$ for
each positive integer $j$.  Thus $T^1 = T$, $T^2 = T \circ T$, etc.,
and it is convenient to take $T^0$ to be the identity mapping $I$ on
$V$.  As in the case of real and complex numbers, one can consider
infinite series of the form
\begin{equation}
        \sum_{j = 0}^\infty T^j,
\end{equation}
and limits of sequences of averages of the form
\begin{equation}
        \frac{I + T + \cdots + T^n}{n}.
\end{equation}
One can also consider these expressions applied to individual vectors
in $V$.

        More precisely, if $V$ has finite dimension, then limits of
sequences of vectors in $V$ can be defined in terms of the
corresponding sequences of coefficients with respect to a basis of
$V$, and limits of linear mappings can be defined in terms of the
entries of the corresponding matrices.  By standard arguments,
convergence is independent of the particular choice of basis of $V$.
Convergence can also be defined with respect to a norm on any vector
space, as in the next section, and is equivalent to using a basis when
$V$ has finite dimension.

        Suppose for instance that $v \in V$ is an eigenvector of $T$
with eigenvalue $\lambda \in {\bf R}$ or ${\bf C}$, as appropriate, so that
\begin{equation}
        T(v) = \lambda \, v.
\end{equation}
For each $j$,
\begin{equation}
        T^j(v) = \lambda^j \, v,
\end{equation}
and we are back to sequences and series of real and complex numbers.
If $V$ has finite dimension and there is a basis of $V$ consisting of
eigenvectors of $T$, then $T^j$ is diagonalized by the same basis for
each $j$, and the previous sequences and series of linear mappings are
reduced to sequences and series of complex numbers.

        Remember that any linear mapping $T$ on a finite-dimensional
complex vector space of positive dimension has a nonzero eigenvector,
as a consequence of the fundamental theorem of algebra.  If $T$ is not
diagonalizable, then the behavior of $T^j$ can still be analyzed in
terms of the Jordan canonical form.

\section{Convergence}
\label{convergence}
\setcounter{equation}{0}

        Let $V$ be a real or complex vector space equipped with a norm
$\|\cdot \|$.  A sequence $\{v_j\}_j$ of vectors in $V$ is said to
\emph{converge} to $v \in V$ if for every $\epsilon > 0$ there is an
$L$ such that
\begin{equation}
        \|v_j - v\| < \epsilon
\end{equation}
for every $j \ge L$.  This is the same as the usual definition of
convergence of a sequence of real or complex numbers when $V = {\bf
R}$ or ${\bf C}$ and the norm is given by the absolute value or
modulus.  As usual, the limit of a sequence is unique when it exists.

        If $\{v_j\}_j$, $\{w_j\}_j$ are sequences of vectors in $V$
that converge to $v, w \in V$, respectively, then the sequence of sums
$v_j + w_j$ converges to $v + w$.  Similarly, if $\{v_j\}_j$ converges
to $v$ in $V$, and $\{t_j\}_j$ is a sequence of real or complex
numbers that converges to $t \in {\bf R}$ or ${\bf C}$, as
appropriate, then $\{t_j \, v_j\}_j$ converges to $t \, v$ in $V$.
These statements can be verified using standard arguments.  An
infinite series $\sum_{j = 0}^\infty a_j$ with terms in $V$ converges
if the corresponding sequence of partial sums $b_n = \sum_{j = 0}^n
a_j$ converges in $V$.  If $\sum_{j = 0}^\infty a_j$ converges, then
$\{a_j\}_j$ converges to $0$ in $V$, as in the case of real or complex
numbers.

        If $V = {\bf R}^n$ or ${\bf C}^n$ with norm $\|v\|_p$ for some
$p$, $1 \le p \le \infty$, then a sequence $\{v_j\}_j$ converges to $v
\in V$ if and only if the corresponding $n$ sequences of coordinates
of the $v_j$'s converge to the coordinates of $v$ as sequences of real
or complex numbers.  This also works for any norm on ${\bf R}^n$ or
${\bf C}^n$, since all norms on these spaces are equivalent.  There
are analogous statements for any finite-dimensional vector space $V$
and any basis in $V$, using a linear mapping that sends the given
basis to the standard basis in ${\bf R}^n$ or ${\bf C}^n$, as
appropriate.

        Suppose that $V$ and $W$ are vector spaces, both real or both
complex, and equipped with norms.  The space of bounded linear
mappings from $V$ into $W$ is also a vector space, and convergence of
sequences and series in this space can be defined in terms of the
operator norm.

\section{Completeness}
\label{completeness}
\setcounter{equation}{0}

        Let $V$ be a real or complex vector space with a norm $\|\cdot
\|$.  A sequence $\{v_j\}_j$ of vectors in $V$ is said to be a
\emph{Cauchy sequence} if for every $\epsilon > 0$ there is an $L$
such that
\begin{equation}
        \|v_j - v_l\| < \epsilon
\end{equation}
for every $j, l \ge L$.  Convergent sequences are Cauchy sequences,
and $V$ is said to be \emph{complete} if every Cauchy sequence in $V$
converges to an element of $V$.  A complete vector space with respect
to a norm is known as a \emph{Banach space}, and it is a \emph{Hilbert
space} if the norm is determined by an inner product.  It is well
known that ${\bf R}$ and ${\bf C}$ are complete with respect to the
usual absolute value and modulus.  A sequence in ${\bf R}^n$ or ${\bf
C}^n$ is a Cauchy sequence with respect to a $p$-norm if and only if
the $n$ sequences of its coordinates are Cauchy sequences in ${\bf R}$
or ${\bf C}$, as appropriate.  It follows that ${\bf R}^n$ and ${\bf
C}^n$ are complete with respect to the $p$-norms, since ${\bf R}$ and
${\bf C}$ are complete.  Hence ${\bf R}^n$ and ${\bf C}^n$ are
complete with respect to any norm, by equivalence of norms.  This
implies in turn that finite-dimensional vector spaces are always
complete.

        Let $\sum_{j = 0}^\infty a_j$ be an infinite series with terms
in $V$.  This series converges \emph{absolutely} if
\begin{equation}
        \sum_{j = 0}^\infty \|a_j\|
\end{equation}
converges as an infinite series of nonnegative real numbers.  The
sequence of partial sums of an absolutely convergent series is a
Cauchy sequence, just as for absolutely convergent series of real or
complex numbers.  If $V$ is a Banach space, then it follows that every
absolutely convergent series of vectors in $V$ converges in $V$.

        Conversely, if every absolutely convergence series in $V$
converges, then $V$ is complete.  For suppose that $\{v_j\}_j$ is a
Cauchy sequence in $V$, and let $\{v_{j_l}\}_{l = 0}^\infty$ be a
subsequence of $\{v_j\}_j$ such that
\begin{equation}
        \|v_{j_l} - v_{j_{l + 1}}\| \le 2^{-l}
\end{equation}
for each $l$.  Thus $\sum_{l = 0}^\infty (v_{j_l} - v_{j_{l + 1}})$
converges absolutely in $V$, and hence converges in $V$ by hypothesis.
This implies that $\{v_{j_l}\}_{l = 0}^\infty$ converges in $V$, since 
\begin{equation}
        \sum_{l = 0}^n (v_{j_l} - v_{j_{l + 1}}) = v_{j_1} - v_{j_{n + 1}}
\end{equation}
for each $n$.  Therefore $\{v_j\}_j$ converges in $V$, because a
Cauchy sequence with a convergent subsequence converges to the same
limit.

        Let $V$ and $W$ be vector spaces, both real or both complex,
and equipped with norms.  If $W$ is complete, then the vector space of
bounded linear mappings from $V$ into $W$ is complete with respect to
the operator norm.  For if $\{T_j\}_j$ is a Cauchy sequence of bounded
linear mappings from $V$ into $W$, then $\{T_j(v)\}_j$ is a Cauchy
sequence in $W$ for each $v \in V$.  Because $W$ is complete,
$\{T_j(v)\}_j$ converges in $W$, and its limit determines a linear
mapping $T$ from $V$ into $W$.  One can check that $T$ is a bounded
linear mapping from $V$ into $W$, and that $\{T_j\}_j$ converges to
$T$ in the operator norm.

\section{The supremum norm}
\label{supremum norm}
\setcounter{equation}{0}

        A continuous real or complex-valued function $f$ on a
topological space $X$ is said to be \emph{bounded} if there is a
nonnegative real number $A$ such that
\begin{equation}
        |f(x)| \le A
\end{equation}
for every $x \in X$.  In this case, the \emph{supremum norm} of $f$ is
defined by
\begin{equation}
        \|f\|_{sup} = \sup \{|f(x)| : x \in X\}.
\end{equation}
This is a norm on the vector space of bounded real or complex-valued
continuous functions on $X$.

        Convergence of a sequence of bounded continuous functions on
$X$ with respect to the supremum norm is the same as uniform
convergence.  If $\{f_j(x)\}_j$ is a Cauchy sequence of bounded
continuous functions on $X$ with respect to the supremum norm, then
$\{f_j(x)\}_j$ is a Cauchy sequence of real or complex numbers for
each $x \in X$, as appropriate.  Hence $\{f_j\}_j$ converges pointwise
to a function $f$ on $X$, and one can use the Cauchy condition with
respect to the supremum norm to show that $\{f_j\}_j$ converges
uniformly to $f$.  This implies that $f$ is bounded and continuous on
$X$, by well-known results about uniform convergence.  Thus
$\{f_j\}_j$ converges to $f$ in the supremum norm, and the space of
bounded continuous functions on $X$ is complete with respect to the
supremum norm.

        Let $V$ be a real or complex vector space equipped with a norm
$\|v\|_V$.  A continuous function $f$ on $X$ with values in $V$ is
said to be \emph{bounded} if $\|f(x)\|_V$ is bounded as a real-valued
function on $X$, in which event the supremum norm of $f$ with respect
to $\|v\|_V$ is defined by
\begin{equation}
        \|f\|_{sup, V} = \sup \{\|f(x)\|_V : x \in X\}.
\end{equation}
This is a norm on the vector space $\mathcal{C}_b(X, V)$ of bounded
continuous functions on $X$ with values in $V$.  If $V$ is complete
with respect to $\|v\|_V$, then $\mathcal{C}_b(X, V)$ is complete with
respect to $\|f\|_{sup, V}$, by an argument like the one in the
previous paragraph.

        For example, if $v \in V$ and $\phi$ is a bounded continuous
real or complex-valued function on $X$, depending on whether $V$ is a
real or complex vector space, then $\Phi(x) = \phi(x) \, v$ is a
bounded continuous function on $X$ with values in $V$, and
\begin{equation}
        \|\Phi\|_{sup, V} = \|\phi\|_{sup} \, \|v\|_V.
\end{equation}
If $f$ is bounded continuous $V$-valued function on $X$, then
$\|f(x)\|_V$ is a bounded continuous real-valued function on $X$ whose
supremum norm is $\|f\|_{sup, V}$.  Of course, every continuous
function on $X$ with values in $V$ is bounded when $X$ is compact.

\section{Algebras}
\label{algebras}
\setcounter{equation}{0}

        Let $\mathcal{A}$ be an associative algebra over the real or
complex numbers.  Thus $\mathcal{A}$ is a vector space over the real
or complex numbers equipped with a binary operation of multiplication
$a \, b$ which is associative in the sense that
\begin{equation}
        (a \, b) \, c = a \, (b \, c)
\end{equation}
for every $a, b, c \in \mathcal{A}$.  More precisely, multiplication
is asked to be a bilinear mapping, which means that $a \mapsto a \, b$
is a linear mapping for each $b \in \mathcal{A}$, and $b \mapsto a \,
b$ is linear for each $a \in \mathcal{A}$.  A pair of elements $a$,
$b$ of $\mathcal{A}$ \emph{commute} with each other if
\begin{equation}
        a \, b = b \, a,
\end{equation}
and $\mathcal{A}$ is a commutative algebra if this holds for every $a,
b \in \mathcal{A}$.

        It will be convenient to suppose also that there be a nonzero
multiplicative identity element $e$ in $\mathcal{A}$, which is to say
that
\begin{equation}
        a \, e = e \, a = a
\end{equation}
for every $a \in \mathcal{A}$.  An element $a$ of $\mathcal{A}$ is
said to be \emph{invertible} if there is a $b \in \mathcal{A}$ such that
\begin{equation}
        a \, b = b \, a = e.
\end{equation}
In this case, the inverse $b$ of $a$ is unique, and denoted $a^{-1}$.
If $a$ is invertible, and $c \in \mathcal{A}$ commutes with $a$, then
$c$ commutes with $a^{-1}$ as well, since
\begin{equation}
        a^{-1} \, c = a^{-1} \, c \, a \, a^{-1}
                    = a^{-1} \, a \, c \, a^{-1}
                    = c \, a^{-1}.
\end{equation}
If $a_1$, $a_2$ are invertible elements of $\mathcal{A}$, then their
product $a_1 \, a_2$ is invertible too, and is given by
\begin{equation}
        (a_1 \, a_2)^{-1} = a_2^{-1} \, a_1^{-1}.
\end{equation}
Conversely, if $a_1$, $a_2$ are commuting elements of $\mathcal{A}$
whose product $a_1 \, a_2$ is invertible, then $a_1$ and $a_2$ are
each invertible.  In this case, $a_1 \, a_2$ and hence its inverse
commute with $a_1$ and $a_2$, and
\begin{equation}
 a_1^{-1} = (a_1 \, a_2)^{-1} \, a_2, \ a_2^{-1} = a_1 \, (a_1 \, a_2)^{-1}.
\end{equation}

        For example, the real or complex-valued continuous functions
on a topological space $X$ form a commutative algebra with respect to
pointwise multiplication of functions.  The constant function equal to
$1$ at every element of $X$ is the multiplicative identity element in
this algebra, and a continuous function $f$ on $X$ is invertible in
this algebra if and only if $f(x) \ne 0$ for every $x \in X$, in which
event the inverse of $f$ is given by $1/f(x)$.  The bounded continuous
function on $X$ also form an algebra which is a subalgebra of the
algebra of all continuous functions on $X$.  In order for a bounded
continuous function $f$ to be invertible in this subalgebra, it is
necessary that there be an $\eta > 0$ such that
\begin{equation}
        |f(x)| \ge \eta
\end{equation}
for every $x \in X$, so that $1/f(x)$ is bounded on $X$.  If $X$ is
compact, then these two algebras are the same.

        If $V$ is a vector space, then the linear mappings on $V$ form
an algebra with composition as multiplication.  The identity operator
$I$ on $V$ is the multiplicative identity element of this algebra, and
a linear mapping $T$ on $V$ is invertible in the algebra if and only
if it is a one-to-one mapping of $V$ onto itself.  If $V$ is equipped
with a norm, then the bounded linear mappings on $V$ with respect to
this norm form an algebra which is a subalgebra of the algebra of all
linear mappings on $V$.  Invertibility of a bounded linear mapping $T$
on $V$ in this subalgebra means that the inverse mapping $T^{-1}$ is
also bounded.  If $V$ has finite dimension, then these two algebras
are the same.

\section{Banach algebras}
\label{banach algebras}
\setcounter{equation}{0}

        Let $\mathcal{A}$ be an associative algebra with nonzero
multiplicative identity element $e$ over the real or complex numbers,
and suppose that $\mathcal{A}$ is equipped with a norm $\|a\|$.  If
\begin{equation}
        \|a \, b\| \le \|a\| \, \|b\|
\end{equation}
for every $a, b \in \mathcal{A}$ and
\begin{equation}
        \|e\| = 1,
\end{equation}
then $(\mathcal{A}, \|a\|)$ is said to be a \emph{normed algebra}.  In
particular, this implies that the product of a pair of convergent
sequences in $\mathcal{A}$ converges to the product of the limits of
the sequences, just as for products of convergent sequences of real or
complex numbers.

        If $\mathcal{A}$ is complete with respect to $\|a\|$, then
$\mathcal{A}$ is said to be a \emph{Banach algebra}.  For example, the
algebra of bounded continuous real or complex-valued functions on a
topological space is a Banach algebra.

        Let $V$ be a vector space over the real or complex numbers
equipped with a norm.  The algebra of bounded linear mappings on $V$
is a normed algebra with respect to the operator norm.  If $V$ is
complete, then the algebra of bounded linear operators on $V$ is a
Banach algebra.

\section{Invertibility}
\label{invertibility}
\setcounter{equation}{0}

        Let $(\mathcal{A}, \|a\|)$ be a normed algebra with nonzero
multiplicative identity element $e$.  For each $a \in \mathcal{A}$ and
positive integer $j$, let $a^j$ be the product of $j$ factors of $a$,
so that $a^1 = a$, $a^2 = a \, a$, etc., with $a^0 = e$.  Thus
\begin{equation}
        \|a^j\| \le \|a\|^j
\end{equation}
for each $j$.  If $\|a\| < 1$, then $\{a^j\}_{j = 0}^\infty$ converges
to $0$ in $\mathcal{A}$,
\begin{equation}
        \sum_{j = 0}^\infty \|a^j\| \le \sum_{j = 0}^\infty \|a\|^j
                                     = \frac{1}{1 - \|a\|},
\end{equation}
and hence $\sum_{j = 0}^\infty a^j$ converges absolutely in $\mathcal{A}$.

        If $\|a\| < 1$ and $\mathcal{A}$ is a Banach algebra, then
it follows that $\sum_{j = 0}^\infty a^j$ converges in $\mathcal{A}$.
For each $n \ge 0$,
\begin{equation}
        (e - a) \Big(\sum_{j = 0}^n a^j\Big)
         = \Big(\sum_{j = 0}^n a^j\Big) (e - a) = e - a^{n + 1},
\end{equation}
which implies that
\begin{equation}
        (e - a) \Big(\sum_{j = 0}^\infty a^j\Big)
         = \Big(\sum_{j = 0}^\infty a^j \Big) (e - a) = e,
\end{equation}
since $a^{n + 1} \to 0$ as $n \to \infty$.  Therefore, $e - a$ is
invertible, with
\begin{equation}
        (e - a)^{-1} = \sum_{j = 0}^\infty a^j.
\end{equation}
This is a fundamental property of Banach algebras.

        Let $x$ be an invertible element of $\mathcal{A}$, so that
\begin{equation}
        1 = \|x \, x^{-1}\| \le \|x\| \, \|x^{-1}\|.
\end{equation}
If $y \in \mathcal{A}$ satisfies
\begin{equation}
        \|x - y\| < \frac{1}{\|x^{-1}\|},
\end{equation}
then
\begin{equation}
        y = x - (x - y) = x \, (e - x^{-1} \, (x - y))
\end{equation}
is invertible by the remarks of the previous paragraph.  In
particular, the invertible elements of $\mathcal{A}$ form an open set.

        If $\|a\| < 1$, then
\begin{equation}
        \|(e - a)^{-1} - e\| = \Bigl\|\sum_{j = 1}^\infty a^j\Bigr\|
          \le \sum_{j = 1}^\infty \|a\|^j = \frac{\|a\|}{1 - \|a\|}.
\end{equation}
Hence $x \mapsto x^{-1}$ is a continuous mapping on the set of
invertible elements of $\mathcal{A}$.

\section{Submultiplicative sequences}
\label{submultiplicative sequences}
\setcounter{equation}{0}

        A sequence $\{r_j\}_{j = 1}^\infty$ of nonnegative real
numbers is said to be \emph{submultiplicative} if
\begin{equation}
\label{submultiplicativity}
        r_{j + l} \le r_j \, r_l
\end{equation}
for every $j, l \ge 1$.  In this case,
\begin{equation}
        \lim_{n \to \infty} r_n^{1/n} = \inf_{n \ge 1} r_n^{1/n},
\end{equation}
where the existence of the limit is part of the conclusion.  To see
this, observe that
\begin{equation}
        r_{j \, n + l} \le (r_n)^j \, (r_1)^l
\end{equation}
and hence
\begin{equation}
        (r_{j \, n + l})^{1/(j \, n + l)}
 \le [(r_n)^{1/n}]^{j \, n / (j \, n + l)} \, (r_1)^{l/(j \, n + l)}
\end{equation}
for every $j, l, n \ge 1$.  If $n$ is fixed, then the right side can
be approximated by $(r_n)^{1/n}$ for $j$ sufficiently large and $0 \le
l < n$.

\section{Invertibility, 2}
\label{invertibility, 2}
\setcounter{equation}{0}

        For each element $x$ of a Banach algebra $\mathcal{A}$, $r_j =
\|x^j\|$ is a submultiplicative sequence of nonnegative real numbers.  Put
\begin{equation}
        \rho(x) = \lim_{n \to \infty} \|x^n\|^{1/n}
                = \inf_{n \ge 1} \|x^n\|^{1/n} \le \|x\|,
\end{equation}
which is known as the \emph{spectral radius} of $x$, at least in the
case of complex Banach algebras.  Thus
\begin{equation}
        \rho(t \, x) = |t| \, \rho(x)
\end{equation}
for each real or complex number $t$, as appropriate.

        If $\rho(x) < 1$, then $\sum_{j = 0}^\infty \|x^j\|$
converges, and hence $\sum_{j = 0}^\infty x^j$ converges in
$\mathcal{A}$.  As before, the sum is equal to the inverse of $e - x$
under these conditions.  Note that $\rho(x) < 1$ if and only if
$\|x^n\|^{1/n} < 1$ for some $n$, which is the same as $\|x^n\| < 1$.
Invertibility of $e - x$ can also be obtained from the invertiblity of
$e - x^n$, since the latter is the product of the former and $e + x +
\cdots + x^{n - 1}$.

        If $\rho(x) \le 1$ and $t$ is a real or complex number,
such that $|t| < 1$, as appropriate, then consider
\begin{equation}
        \sum_{j = 0}^\infty t^j \, x^j = (e - t \, x)^{-1}.
\end{equation}
If $e - x$ is invertible, then $(e - t \, x)^{-1}$ extends
continuously to a neighborhood of $t = 1$.  Conversely, if there is a
sequence $\{t_j\}_j$ such that $|t_j| < 1$ for each $j$, $\{t_j\}_j$
converges to $1$, and $(e - t_j \, x)^{-1}$ converges in
$\mathcal{A}$, then $e - x$ is invertible and the limit of $(e - t_j
\, x)^{-1}$ is equal to its inverse.

        Suppose that $\mathcal{A}$ is the algebra of bounded
continuous real or complex-valued functions on a topological space
$X$, with the supremum norm.  In this case,
\begin{equation}
        \|f^n\|_{sup} = \|f\|_{sup}^n
\end{equation}
for each $n \ge 1$ and $f$, and hence $\rho(f) = \|f\|_{sup}$.

\section{Spectrum}
\label{spectrum}
\setcounter{equation}{0}

        Let $T$ be a linear mapping on a real or complex vector space
$V$.  To say that a real or complex number $\lambda$, as appropriate,
is an eigenvalue of $T$ means exactly that $T - \lambda \, I$ has
nontrivial kernel.  In particular, $T - \lambda \, I$ is not
invertible.  Conversely, if $V$ has finite dimension and $T - \lambda
\, I$ is not invertible, then $T - \lambda \, I$ has nontrivial
kernel, and $\lambda$ is an eigenvalue of $T$.

        Suppose that $V$ is equipped with a norm $\|\cdot \|_V$, and
that $T$ is a bounded linear operator on $V$.  If $v$ is a nonzero
eigenvector of $T$ with eigenvalue $\lambda$, then
\begin{equation}
        |\lambda| \, \|v\|_V = \|T(v)\|_V \le \|T\|_{op} \, \|v\|_V
\end{equation}
implies that
\begin{equation}
        |\lambda| \le \|T\|_{op}.
\end{equation}
Similarly, $\lambda^n$ is an eigenvalue of $T^n$ for each positive
integer $n$, and hence
\begin{equation}
        |\lambda|^n \le \|T^n\|_{op}.
\end{equation}
Thus
\begin{equation}
        |\lambda| \le \|T^n\|_{op}^{1/n}
\end{equation}
for each $n$, and therefore $|\lambda| \le \rho(T)$.

        Let $(\mathcal{A}, \|\cdot \|)$ be a Banach algebra with
nonzero multiplicative identity element $e$, let $x$ be an element of
$\mathcal{A}$, and let $\lambda$ be a real or complex number, as
appropriate.  If $|\lambda| > \|x\|$, then
\begin{equation}
        \lambda \, e - x = \lambda \, (e - \lambda^{-1} \, x)
\end{equation}
is invertible.  The same conclusion holds when $|\lambda| > \rho(x)$.
Equivalently,
\begin{equation}
        |\lambda| \le \rho(x)
\end{equation}
when $\lambda \, e - x$ is not invertible in $\mathcal{A}$.

        The set
\begin{equation}
        \sigma(x) = \{\lambda : \lambda \, e - x \hbox{ is not invertible in }
                                                   \mathcal{A}\}
\end{equation}
is known as the \emph{spectrum} of $x$ in $\mathcal{A}$, especially in
the complex case.  This is a closed set in ${\bf R}$ or ${\bf C}$, as
appropriate, since the complementary \emph{resolvent} set of $\lambda$
such that $\lambda \, e - x$ is invertible is an open set when
$\mathcal{A}$ is a Banach algebra.  If $\mathcal{A}$ is a complex
Banach algebra, then a famous theorem states that $\sigma(x) \ne
\emptyset$ for every $x \in \mathcal{A}$.  Basically, if $\sigma(x) =
\emptyset$, then
\begin{equation}
        (\lambda \, e - x)^{-1}
\end{equation}
would be a nonconstant holomorphic $\mathcal{A}$-valued function on
the complex plane that tends to $0$ as $|\lambda| \to \infty$, a contradiction.
This is still a holomorphic $\mathcal{A}$-valued function on the complement
of $\sigma(x)$ in the complex plane for any $x \in \mathcal{A}$.
Another famous theorem states that
\begin{equation}
        \rho(x) = \max \{|\lambda| : \lambda \in \sigma(x)\}
\end{equation}
when $\mathcal{A}$ is a complex Banach algebra, as a consequence of the
convergence of
\begin{equation}
        \sum_{j = 0}^\infty \alpha^j \, x^j
\end{equation}
when $\alpha$ is a nonzero complex number such that $|\lambda| <
1/|\alpha|$ for every $\lambda \in \sigma(x)$.

\section{Averages in normed algebras}
\label{averages in normed algebras}
\setcounter{equation}{0}

        Let $(\mathcal{A}, \|\cdot \|)$ be a normed algebra with
nonzero multiplicative identity element $e$.  If $x \in \mathcal{A}$
and $\|x\| < 1$, then
\begin{equation}
        \Bigl\|\sum_{j = 0}^n x^j\Bigr\| \le \sum_{j = 0}^n \|x^j\|
         \le \sum_{j = 0}^n \|x\|^j \le \frac{1}{1 - \|x\|}
\end{equation}
for each $n$.  Hence
\begin{equation}
\label{lim (1/n+1) sum_{j = 0}^n x^j = 0}
        \lim_{n \to\infty} \frac{1}{n + 1} \sum_{j = 0}^n x^j = 0.
\end{equation}

        Suppose now that $\|x\| = 1$.  Thus
\begin{equation}
        \Bigl\|\frac{1}{n + 1} \sum_{j = 0}^n x^j\Bigr\|
         \le \frac{1}{n + 1} \sum_{j = 0}^n \|x^j\|
         \le \frac{1}{n + 1} \sum_{j = 0}^n \|x\|^j = 1
\end{equation}
for each $n$.  Of course,
\begin{equation}
        \frac{1}{n + 1} \sum_{j = 0}^n x^j = e
\end{equation}
for each $n$ when $x = e$.

        For any $x \in \mathcal{A}$ and $n \ge 1$,
\begin{equation}
        (e - x) \sum_{j = 0}^n x^j = e - x^{n + 1}.
\end{equation}
This implies that
\begin{equation}
        \sum_{j = 0}^n x^j = (e - x)^{-1} \, (e - x^{n + 1})
\end{equation}
when $e - x$ is invertible in $\mathcal{A}$.  If $\|x\| = 1$ and $e -
x$ is invertible, then
\begin{equation}
        \Bigl\|\sum_{j = 0}^n x^j\Bigr\|
          \le \|(e - x)^{-1}\| \, \|e - x^{n + 1}\|
           \le 2 \, \|(e - x)^{-1}\|,
\end{equation}
so that (\ref{lim (1/n+1) sum_{j = 0}^n x^j = 0}) holds in this case too.

        If $e - x$ is invertible, then
\begin{eqnarray}
        \sum_{l = 0}^n \sum_{j = 0}^l x^j 
         & = & \sum_{l = 0}^n (e - x)^{-1} \, (e - x^{l + 1}) \\
 & = & (n + 1) \, (e - x)^{-1} - (e - x)^{-2} \, (x - x^{n + 2}) \nonumber
\end{eqnarray}
for each $n$.  If in addition $\|x\| = 1$, then it follows that
\begin{equation}
        \lim_{n \to \infty} \frac{1}{n + 1} \sum_{l = 0}^n \sum_{j = 0}^l x^j
         = (e - x)^{-1}.
\end{equation}

\section{Invertibility, 3}
\label{invertibility, 3}
\setcounter{equation}{0}

        Let $\mathcal{A}$ be a Banach algebra with nonzero
multiplicative identity element $e$.  Suppose that $x \in \mathcal{A}$
has the property that the sums
\begin{equation}
\label{sum_{j = 0}^n x^j}
        \sum_{j = 0}^n x^j
\end{equation}
are uniformly bounded in $\mathcal{A}$.  In particular, $x^0 = e, \,
x, x^2, \ldots$ is a bounded sequence in $\mathcal{A}$, so that
\begin{equation}
\label{sum_{j = 0}^infty r^j x^j}
        \sum_{j = 0}^\infty r^j \, x^j
\end{equation}
converges absolutely when $0 \le r < 1$.  Using summation by parts,
we get that (\ref{sum_{j = 0}^infty r^j x^j}) is equal to
\begin{equation}
        (1 - r) \sum_{l = 0}^\infty r^l \Big(\sum_{j = 0}^l x^j\Big).
\end{equation}
Thus the boundedness of (\ref{sum_{j = 0}^n x^j}) implies the
boundedness of (\ref{sum_{j = 0}^infty r^j x^j}).  Hence the inverse
of $e - r \, x$ has bounded norm when $0 \le r < 1$, since it is given
by (\ref{sum_{j = 0}^infty r^j x^j}).  It follows that $e - x$ is
invertible, by the results of Section \ref{invertibility}.

\section{The open mapping theorem}
\label{open mapping theorem}
\setcounter{equation}{0}

        Let $(V, \|\cdot \|_V)$, $(W, \|\cdot \|_W)$ be Banach spaces,
both real or both complex, and let $T$ be a bounded linear mapping
from $V$ into $W$.  If $T$ maps $V$ \emph{onto} $W$, then Banach's
open mapping theorem says that $T$ sends open subsets of $V$ to open
subsets of $W$.  In particular, if $T$ is a one-to-one mapping of $V$
onto $W$, then it follows that $T^{-1}$ is a bounded linear mapping
from $W$ onto $V$.

        For each $r > 0$, put
\begin{equation}
        B_V(r) = \{v \in V : \|v\|_V < r\},
\end{equation}
and let $B_W(r)$ be defined in the same way.  Because $T$ is linear,
it suffices to show that there is an $r > 0$ such that
\begin{equation}
        T(B_V(1)) \supseteq B_W(r).
\end{equation}

        The hypothesis that $T$ map $V$ onto $W$ implies that
\begin{equation}
        \bigcup_{n = 1}^\infty T(B_V(n)) = W.
\end{equation}
Hence $W$ is the union of the closure $\overline{T(B_V(n))}$ of
$T(B_V(n))$, $n \ge 1$, and it follows from the Baire category theorem
that $\overline{T(B_V(n))}$ contains a nonempty open set for some $n$.
Using linearity, one can show that there is an $r_1 > 0$ such that
\begin{equation}
\label{overline{T(B_V(1))} supseteq B_W(r_1)}
        \overline{T(B_V(1))} \supseteq B_W(r_1).
\end{equation}

        Let us use completeness of $V$ and linearity of $T$ to show that
\begin{equation}
\label{T(B_V(2)) supseteq B_W(r_1)}
        T(B_V(2)) \supseteq B_W(r_1).
\end{equation}
Let $w \in B_W(r_1)$ be given.  By (\ref{overline{T(B_V(1))} supseteq
B_W(r_1)}), there is a $v_0 \in B_V(1)$ such that
\begin{equation}
        \|w - T(v_0)\|_W < \frac{r_1}{2}.
\end{equation}
Applying the same argument to $2 (w - T(v_0))$, we get that there
is a $v_1 \in B_V(1/2)$ such that
\begin{equation}
        \|w - T(v_0) - T(v_1)\|_W < \frac{r_1}{4}.
\end{equation}
Repeating the process, we get $v_0, v_1, v_2, \ldots \in V$ such that
\begin{equation}
        \|v_j\|_V < 2^{-j}
\end{equation}
for each $j$, and
\begin{equation}
        \Bigl\|w - \sum_{j = 0}^l T(v_j)\Bigr\| < 2^{-l-1} \, r_1
\end{equation}
for each $l$.  Thus
\begin{equation}
        \sum_{j = 0}^\infty \|v_j\|_V < 2,
\end{equation}
and so $\sum_{j = 0}^\infty v_j$ converges in $V$ and satisfies
\begin{equation}
        \Bigl\|\sum_{j = 0}^\infty v_j\Bigr\|_V < 2.
\end{equation}
Moreover,
\begin{equation}
        T\Big(\sum_{j = 0}^\infty v_j\Big) = w,
\end{equation}
which implies (\ref{T(B_V(2)) supseteq B_W(r_1)}).

\section{The uniform boundedness principle}
\label{uniform boundedness principle}
\setcounter{equation}{0}

        Let $V$ and $W$ be vector spaces, both real or both complex,
equipped with norms $\|\cdot \|_V$ and $\|\cdot \|_W$, respectively.
Let $T_1, T_2, \ldots$ be a sequence of bounded linear mappings from
$V$ into $W$ such that $\{T_j(v)\}_{j = 1}^\infty$ is a bounded
sequence in $W$ for each $v \in V$.  If $V$ is complete, then a
theorem of Banach and Steinhaus implies that the operator norms of the
$T_j$'s are bounded.

        As a variant of this, let $f_1, f_2, \ldots$ be a sequence of
nonnegative real-valued continuous functions on a metric space $M$
such that $\{f_j(x)\}_{j = 1}^\infty$ is bounded for each $x \in M$.
If $M$ is complete, then there is a nonempty open set in $M$ on which
the $f_j$'s are uniformly bounded.  To see this, consider
\begin{equation}
        E_n = \{x \in M : f_j(x) \le n \hbox{ for each } j\},
\end{equation}
which is a closed set in $M$ for each $n$ by continuity.  The
hypothesis of pointwise boundedness means exactly that
\begin{equation}
        \bigcup_{n = 1}^\infty E_n = M,
\end{equation}
and the Baire category theorem implies that $E_n$ contains a nonempty
open set for some $n$, as desired.

        Let us apply this to $M = V$ and $f_j(v) = \|T_j(v)\|_W$.
Note that $\|T_j(v)\|_W$ is a cintinuous function on $V$, since $T_j :
V \to W$ is a bounded linear mapping.  Because of linearity, uniform
boundedness of $\|T_j(v)\|_W$ on a nonempty open set in $V$ implies
that the operator norms of the $T_j$'s are bounded.

        In practice, we are especially interested in situations where
$\{T_j(v)\}_{j = 1}^\infty$ converges in $W$ for each $v \in V$.

\section{Strong operator convergence}
\label{strong operator convergence}
\setcounter{equation}{0}

        Let $V$ and $W$ be vector spaces, both real or both complex,
with norms $\|\cdot \|_V$ and $\|\cdot \|_W$.  A sequence $T_1, T_2,
\ldots$ of bounded linear mappings from $V$ into $W$ is said to
converge \emph{strongly} to a linear mapping $T$ from $V$ into $W$ if
the following two conditions are satisfied.  First, the operator norms
of the $T_j$'s are uniformly bounded, so that there is an $A \ge 0$
such that
\begin{equation}
        \|T_j\|_{op} \le A
\end{equation}
for each $j$.  Second,
\begin{equation}
        \lim_{j \to \infty} T_j(v) = T(v)
\end{equation}
in $W$ for each $v \in V$.  It follows that $T$ is a bounded linear
mapping from $V$ into $W$, with
\begin{equation}
        \|T\|_{op} \le A.
\end{equation}
Of course, convergence in the operator norm implies strong
convergence.  If $V$ is complete, then the pointwise convergence of a
sequence of bounded linear mappings on $V$ implies the boundedness of
their operator norms, as in the previous section.  Conversely, if the
operator norms of the $T_j$'s are uniformly bounded, and if $T$ is a
bounded linear mapping from $V$ into $W$, then convergence of $T_j(v)$
to $T(v)$ for every $v$ in a dense set in $V$ implies the same
property for every $v \in V$.  Similarly, if the operator norms of the
$T_j$'s are uniformly bounded, and if $\{T_j(v)\}_{j = 1}^\infty$ is a
Cauchy sequence in $W$ for each $v$ in a dense set in $V$, then
$\{T_j(v)\}_{j = 1}^\infty$ is a Cauchy sequence in $W$ for every $v
\in V$.  If $W$ is complete, then it follows that $\{T_j\}_{j =
1}^\infty$ converges strongly to a bounded linear mapping from $V$
into $W$.

\section{Convergence of averages}
\label{convergence of averages}
\setcounter{equation}{0}

        Let $V$ be a vector space with a norm $\|\cdot \|$, and let
$T$ be a bounded linear operator on $V$ with $\|T\|_{op} \le 1$.
Thus
\begin{equation}
        \Bigl\|\frac{1}{n + 1} \sum_{j = 0}^n T^j\Bigr\|_{op}
         \le \frac{1}{n + 1} \sum_{j = 0}^n \|T^j\|_{op}
          \le \frac{1}{n + 1} \sum_{j = 0}^n \|T\|_{op}^j \le 1
\end{equation}
for each $n$.

        For $v \in V$, let us consider the convergence in $V$ of
\begin{equation}
\label{1/n + 1 sum_{j = 0}^n T^j(v)}
        \frac{1}{n + 1} \sum_{j = 0}^n T^j(v).
\end{equation}
If $T(v) = v$, then $T^j(v) = v$ for each $j$, and (\ref{1/n + 1
sum_{j = 0}^n T^j(v)}) is equal to $v$ for every $n$.

        If $v = T(u) - u$ for some $u \in V$, then
\begin{equation}
        \sum_{j = 0}^n T^j(v) = T^{n + 1}(u) - u
\end{equation}
for each $n$.  Hence (\ref{1/n + 1 sum_{j = 0}^n T^j(v)}) tends to $0$
as $n \to \infty$, since $T^{n + 1}(u)$ is bounded.  The same conclusion
holds when $v$ is in the closure of the set of $T(u) - u$, $u \in V$.

        Observe that
\begin{equation}
        T\Big(\frac{1}{n + 1}\sum_{j = 0}^n T^j(v)\Big)
              - \frac{1}{n + 1} \sum_{j = 0}^n T^j(v)
        = \frac{1}{n + 1}(T^{n + 1}(v) - v)
\end{equation}
converges to $0$ as $n \to \infty$ for every $v \in V$.  If (\ref{1/n
+ 1 sum_{j = 0}^n T^j(v)}) converges for some $v \in V$, then it
follows that the limit is an eigenvector of $T$ with eigenvalue $1$.

\section{Hilbert spaces}
\label{hilbert spaces}
\setcounter{equation}{0}

        Let $(V, \langle \cdot, \cdot \rangle)$ be a Hilbert space,
and let $E$ be a nonempty closed convex set in $V$.  Let $v \in V$ be
given, and consider
\begin{equation}
        r = \inf \{\|v - w\| : w \in E\}.
\end{equation}
Let $w_1, w_2, \ldots$ be a sequence of elements of $E$ such that
\begin{equation}
        \lim_{j \to \infty} \|v - w_j\| = r.
\end{equation}
Because $E$ is convex, $(w_j + w_l)/2 \in E$ for every $j$, $l$,
and hence
\begin{equation}
        \Bigl\|v - \frac{w_j + w_l}{2}\Bigr\| \ge r.
\end{equation}
The parallelogram law
\begin{equation}
 \|a - b\|^2 + \|a + b\|^2 = 2 \, \|a\|^2 + 2 \, \|b\|^2
\end{equation}
with $a = v - w_j$, $b = v + w_j$ implies that
\begin{equation}
        \|w_j - w_l\|^2 + \Bigl\|2 \, v - w_j - w_l\|^2
         = 2 \, \|v - w_j\|^2 + 2 \, \|v - w_l\|^2.
\end{equation}
It follows that
\begin{equation}
        \lim_{j, l \to \infty} \|w_j - w_l\| = 0,
\end{equation}
which is to say that $\{w_j\}_{j = 1}^\infty$ is a Cauchy sequence in $V$.
This sequence converges, because $V$ is complete, and its limit $w$
is an element of $E$ and satisfies
\begin{equation}
        \|v - w\| = r.
\end{equation}

        Let us apply this to a closed linear subspace $W$ of $V$.  For
each $v \in V$, the preceding argument implies that there is a $w \in
W$ such that
\begin{equation}
        \|v - w\| \le \|v - w - z\|
\end{equation}
for every $z \in W$.  By standard computations, this implies in turn
that
\begin{equation}
        \langle v - w, z \rangle = 0
\end{equation}
for every $z \in W$, in the same way that the minimum of a function is
attained at a critical point.  Conversely, the latter condition implies that
\begin{equation}
        \|v - w - z\|^2 = \|v - w\|^2 + \|z\|^2 \ge \|v - w\|^2
\end{equation}
for every $z \in W$.  If $u \in W$ also satisfies $\langle v - u, z
\rangle = 0$ for every $z \in W$, then
\begin{equation}
        \quad \|u - w\|^2 = \langle u - w, u - w \rangle
 = \langle u - v, u - w\rangle + \langle v - w, u - w \rangle = 0,
\end{equation}
since $u - w \in W$, and hence $u = w$.

        Consider the closed linear subspace $W^\perp$ of $V$ defined by
\begin{equation}
        W^\perp = \{y \in V : \langle y, z \rangle = 0
                               \hbox{ for every } z \in W\}.
\end{equation}
The previous arguments show that every element of $V$ can be expressed
in a unique way as the sum of an element of $W$ and an element of
$W^\perp$.  Let us check that
\begin{equation}
        W = (W^\perp)^\perp.
\end{equation}
Every element of $W$ is contained in $(W^\perp)^\perp$ by definition,
and so it suffices to check that $x \in (W^\perp)^\perp$ is in $W$.
If $x = w + y$ for some $w \in W \subseteq (W^\perp)^\perp$ and $y \in
W^\perp$, then it follows that $y \in (W^\perp)^\perp$, so that
$\langle y, y \rangle = 0$, $y = 0$, and $x = w$.  If $W$ is any
linear subspace of $V$, then $W^\perp$ can be defined in the same way,
and is a closed linear subspace of $V$.  Of course, the closure $\overline{W}$
of $W$ is a closed linear subspace of $V$, and one can check that
\begin{equation}
        \overline{W}^\perp = W^\perp.
\end{equation}
Therefore $\overline{W} = (\overline{W}^\perp)^\perp =
(W^\perp)^\perp$.

        If $W$ is a closed linear subspaces of $V$, then the
\emph{orthogonal projection} $P = P_W$ of $V$ onto $W$ is the linear
mapping on $V$ defined by
\begin{equation}
        P(v) = w
\end{equation}
for $v = w + y$ with $w \in W$, $y \in W^\perp$.  Thus
\begin{equation}
        \|P(v)\| \le \|v\|
\end{equation}
for every $v \in V$, since
\begin{equation}
        \|v\|^2 = \|w\|^2 + \|y\|^2
\end{equation}
in this case.  In particular, $\|P\|_{op} = 1$ except when $W = \{0\}$
and $P = 0$.

\section{Unitary transformations}
\label{unitary transformations}
\setcounter{equation}{0}

        Let $(V, \langle \cdot, \cdot \rangle)$ be an inner product
space.  A one-to-one linear mapping $T$ from $V$ onto itself is said
to be \emph{unitary} if
\begin{equation}
\label{langle T(v), T(w) rangle = langle v, w rangle}
        \langle T(v), T(w) \rangle = \langle v, w \rangle
\end{equation}
for every $v, w \in V$.  This may also be described as an orthogonal
transformation in the real case.  If $T$ is unitary, then
\begin{equation}
        \|T(v)\| = \|v\|
\end{equation}
for every $v \in V$, by taking $v = w$ in the previous equation.
Conversely, a linear mapping $T$ of $V$ onto itself which preserves
norms is one-to-one and satisfies (\ref{langle T(v), T(w) rangle =
langle v, w rangle}), by a polarization argument.  Unitary
transformations are obviously bounded, with operator norm $1$.  Note
that the inverse of a unitary transformation is also unitary.

        Suppose now that $V$ is complete, which is to say that $V$ is
a Hilbert space.  Consider the linear subspace
\begin{equation}
        W = \{T(u) - u : u \in V\}.
\end{equation}
Thus $T(W) = W$, since $T(V) = V$.  Because $T$ is unitary,
\begin{equation}
        \langle y, T(u) - u \rangle = \langle T^{-1}(y) - y, u \rangle
\end{equation}
for every $u, y \in V$.  It follows that $W^\perp$ consists of the $y
\in V$ such that $T^{-1}(y) = y$, which is equivalent to $T(y) = y$.

        As in the previous section, every $v \in V$ can be expressed
in a unique way as $w + y$, where $w \in \overline{W}$, $y \in
W^\perp$, and hence $T(y) = y$.  This implies that
\begin{equation}
        \frac{1}{n + 1} \sum_{j = 0}^n T^j(v)
         = \frac{1}{n + 1} \sum_{j = 0}^n T^n(w) + y
\end{equation}
converges to $y$ in $V$ as $n \to \infty$, as in Section
\ref{convergence of averages}.

        Let $T_W$ be the restriction of $T$ to $\overline{W}$, and let
$I_W$ be the identity mapping on $\overline{W}$.  If $I_W - T_W$ has a
bounded inverse on $\overline{W}$, then
\begin{equation}
        \frac{1}{n + 1} \sum_{j = 0}^n (T_W)^j
\end{equation}
converges to $0$ in the operator norm on $\overline{W}$ as $n \to
\infty$, as in Section \ref{averages in normed algebras}.  In this
case,
\begin{equation}
        \frac{1}{n + 1} \sum_{j = 0}^n T^j
\end{equation}
converges to the orthogonal projection of $V$ onto $W^\perp$ in the
operator norm.

\section{Rotations}
\label{rotations}
\setcounter{equation}{0}

        Put
\begin{equation}
        {\bf T} = \{\alpha \in {\bf C} : |\alpha| = 1\}.
\end{equation}
For each $\alpha \in {\bf T}$, let $R_\alpha$ be the linear mapping on
the vector space of continuous complex-valued functions on ${\bf T}$
defined by
\begin{equation}
        R_\alpha(f)(z) = f(\alpha \, z)
\end{equation}
for each $z \in {\bf T}$.  This is an isometry of the space of
continuous functions on ${\bf T}$ onto itself with respect to the
supremum norm.

        If $f(z) = z^l$ for some integer $l$, then
\begin{equation}
        R_\alpha(f) = \alpha^l \, f,
\end{equation}
so that $f$ is an eigenvector of $R_\alpha$ with eigenvalue $\alpha^l$.
Note that $f(z)$ can also be expressed as $\overline{z}^{-l}$,
since $|z| = 1$.  It follows that
\begin{equation}
\label{1/n+1 sum_{j = 0}^n R_alpha^j(f)}
        \frac{1}{n + 1} \sum_{j = 0}^n R_\alpha^j(f)
\end{equation}
converges uniformly on ${\bf T}$ as $n \to \infty$ when $f$ is one of these
eigenfunctions, or a finite linear combination of these eigenfunctions.

        It is well known that the finite linear combinations of the
functions $z^l$, $l \in {\bf Z}$, are dense in the space of continuous
functions on ${\bf T}$ with respect to the supremum norm.  This
implies that (\ref{1/n+1 sum_{j = 0}^n R_alpha^j(f)}) converges
uniformly on ${\bf T}$ for every continuous function $f$ on ${\bf T}$.

        One can also consider the Lebesgue spaces $L^p({\bf T})$ of
measurable complex-valued functions $f$ on ${\bf T}$ such that $|f|^p$
is integrable with respect to Lebesgue measure on ${\bf T}$, $1 \le p
< \infty$.  For each $\alpha \in {\bf T}$, $R_\alpha$ determines an
isometric linear mapping of $L^p({\bf T})$ onto itself, which is a
unitary mapping when $p = 2$.  The sum (\ref{1/n+1 sum_{j = 0}^n
R_alpha^j(f)}) also converges in the $L^p$ norm for every $f \in
L^p({\bf T})$ when $1 \le p < \infty$, since it converges for a dense
class of functions.

\section{Fourier series}
\label{fourier series}
\setcounter{equation}{0}

        If $f$ is a continuous complex-valued function on the unit
circle ${\bf T}$ and $\alpha \in {\bf T}$, then
\begin{equation}
        \int_{\bf T} f(\alpha \, z) \, |dz| = \int_{\bf T} f(z) \, |dz|.
\end{equation}
Here $|dz|$ denotes the element of arc length integration on ${\bf
T}$, which is invariant under rotations.  This implies that
\begin{equation}
        \int_{\bf T} z^l \, |dz| = 0
\end{equation}
for every nonzero integer $l$.

        For any pair of continuous functions $f_1$, $f_2$ on ${\bf T}$, put
\begin{equation}
        \langle f_1, f_2 \rangle_{\bf T}
         = \frac{1}{2 \pi} \int_{\bf T} f_1(z) \, \overline{f_2(z)} \, |dz|.
\end{equation}
This defines an inner product on the vector space of complex-valued
continuous functions on ${\bf T}$, which extends to the Lebesgue space
$L^2({\bf T})$ of square-integrable functions on ${\bf T}$, so that
the latter becomes a Hilbert space.

        Observe that
\begin{eqnarray}
 \langle z^j, z^l \rangle_{\bf T}  & = & 0 \quad\hbox{when } j \ne l \\
                                   & = & 1 \quad\hbox{when } j = l. \nonumber
\end{eqnarray}
Thus the functions $z^l$, $l \in {\bf Z}$, are orthonormal with respect
to this inner product.

        The Fourier coefficients of a function $f$ on ${\bf T}$ are defined by
\begin{equation}
        \widehat{f}(j) = \langle f, z^j \rangle_{\bf T}
                       = \frac{1}{2 \pi} \int_{\bf T} f(z) \, z^{-j} \, |dz|,
\end{equation}
and the Fourier series associated to $f$ is
\begin{equation}
        \sum_{j = -\infty}^\infty \widehat{f}(j) \, z^j.
\end{equation}
It is well known that the Abel sums of the Fourier series of a
continuous function $f$ on ${\bf T}$, given by
\begin{equation}
 \sum_{j = -\infty}^\infty \widehat{f}(j) \, r^{|j|} \, z^j, \quad 0 < r < 1,
\end{equation}
converge to $f(z)$ uniformly on ${\bf T}$ as $r \to 1$, because these sums
can be expressed as averages of $f$ concentrated near $z$ on ${\bf T}$.
Similarly, the Cesaro means
\begin{equation}
 \frac{1}{n + 1} \sum_{l = 0}^n \sum_{j = -l}^l \widehat{f}(j) \, z^j
\end{equation}
converge to $f(z)$ uniformly on ${\bf T}$ as $n \to \infty$.

\section{Measure-preserving transformations}
\label{meaure-preserving transformations}
\setcounter{equation}{0}

        Let $X$ be a set with a $\sigma$-algebra of measurable subsets
and a positive measure $\mu$.  Let $\phi$ be a one-to-one mapping of
$X$ onto itself such that $E \subseteq X$ is measurable if and only if
$\phi(E)$ is measurable, and
\begin{equation}
        \mu(\phi(E)) = \mu(E)
\end{equation}
for every measurable set $E \subseteq X$.  Conside the linear mapping
\begin{equation}
        T(f) = f \circ \phi
\end{equation}
acting on measurable functions $f$ on $X$.  If $f$ is an integrable
function on $X$, then $T(f)$ is too, and
\begin{equation}
        \int_X T(f) \, d\mu = \int_X f \, d\mu.
\end{equation}
Similarly, $T$ determines an isometry of $L^p(X)$ onto itself for each
$p$, $1 \le p \le \infty$.

        In particular, $T$ defines a unitary mapping on the Hilbert
space $L^2(X)$.  As in Section \ref{unitary transformations},
\begin{equation}
        \frac{1}{n + 1} \sum_{j = 0}^n T^j(f)
\end{equation}
converges in $L^2(X)$ for each $f \in L^2(X)$.  For any $p$,
\begin{equation}
        \Bigl\|\frac{1}{n + 1} \sum_{j = 0}^n T^j(f)\Bigr\|_p
         \le \frac{1}{n + 1} \sum_{j = 0}^n \|T^j(f)\|_p = \|f\|_p
\end{equation}
for each $f \in L^p(X)$ and $n \ge 0$, since $T$ is an isometry on
$L^p(X)$.  If $f$ is in $L^1(X)$ and $L^\infty(X)$, then $f \in
L^p(X)$ for $1 < p < \infty$, and the previous statements imply that
the averages converge in $L^p(X)$ for $1 < p < \infty$.  Because
$L^1(X) \cap L^\infty(X)$ is dense in $L^p(X)$ for $1 < p < \infty$,
it follows that the averages converge in $L^p(X)$ for every $f \in
L^p(X)$ when $1 < p < \infty$.

        If $\mu(X) < \infty$, then $L^2(X)$ is already a dense linear
subspace of $L^1(X)$.  In this case, convergence in $L^2(X)$ implies
convergence in $L^1(X)$, and so the averages converge in $L^1(X)$ for
every $f \in L^2(X)$.  Hence the averages converge in $L^1(X)$ for
every $f \in L^1(X)$.

        Suppose that $X$ is the set ${\bf Z}$ of integers, and that
$\mu$ is counting measure on ${\bf Z}$.  If $\phi(l) = l + 1$ for each
$n$ and $f$ is a function equal to $1$ at one point and to $0$
elsewhere, then the averages have $L^1$ norm equal to $1$ for each
$n$, but converge to $0$ uniformly on $X$.

\section{Sequence spaces}
\label{sequence spaces}
\setcounter{equation}{0}

        Let $A$ be a finite set with at least two elements, and let
$X$ be the set of doubly-infinite sequences $x = \{x_j\}_{j =
-\infty}^\infty$ with $x_j \in A$ for each $j$.  Thus $X$ is a compact
Hausdorff space with respect to the product topology defined using the
discrete topology on $A$.  Let $\phi : X \to X$ be the shift mapping
such that the $j$th term of $\phi(x)$ is equal to $x_{j + 1}$.  This
is a homeomorphism from $X$ onto itself.  Suppose that to each $a \in
A$ is assigned a nonnegative real number $w(a)$ such that
\begin{equation}
        \sum_{a \in A} w(a) = 1.
\end{equation}
This defines a probability measure on $A$, which leads to a
probability measure on $X$.  This probability measure on $X$ is
preserved by the shift mapping $\phi$.

\section{Maximal functions}
\label{maximal functions}
\setcounter{equation}{0}

        Let $f$ be a locally-integrable function on the real line.
Consider the Hardy--Littlewood maximal function $f^*$ associated to
$f$, defined by
\begin{equation}
        f^*(x) = \sup_{x \in I} \frac{1}{|I|} \int_I |f(y)| \, dy,
\end{equation}
where the supremum is taken over all open intervals $I = (a, b)$ in ${\bf R}$
that contain $x$, and where $|I|$ denotes the length $b - a$ of $I$.
Thus $f^*(x) \ge 0$ by construction, and $f^*(x) = +\infty$ is certainly
possible.  The mapping $f \mapsto f^*$ is sublinear in the sense that
\begin{equation}
        (f_1 + f_2)^* \le f_1^* + f_2^*
\end{equation}
for any locally-integrable functions $f$ on ${\bf R}$, and 
\begin{equation}
        (t \,f)^* = |t| \, f^*
\end{equation}
for any locally-integrable function $f$ on ${\bf R}$ and constant $t$.
If $f$ is bounded, then $f^*$ is bounded, and
\begin{equation}
\label{||f^*||_infty le ||f||_infty}
        \|f^*\|_\infty \le \|f\|_\infty.
\end{equation}

        If $f$ is any locally-integrable function on ${\bf R}$, then
$f^*$ is lower semicontinuous.  This means that
\begin{equation}
        E_\lambda = \{x \in {\bf R} : f^*(x) > \lambda\}
\end{equation}
is an open set for each $\lambda \ge 0$.  For if $f^*(x) > \lambda$
for some $x \in {\bf R}$ and $\lambda \ge 0$, then
\begin{equation}
\label{frac{1}{|I|} int_I |f(y)| dy > lambda}
        \frac{1}{|I|} \int_I |f(y)| \, dy > \lambda
\end{equation}
for some open interval $I$ in ${\bf R}$ that contains $x$.  This
implies in turn that $f^* > \lambda$ on $I$, amd hence that $I
\subseteq E_\lambda$, as desired.

        If $f$ is an integrable function on ${\bf R}$, then
\begin{equation}
\label{|E_lambda| le frac{2}{lambda} int_{bf R} |f(y)| dy}
        |E_\lambda| \le \frac{2}{\lambda} \int_{\bf R} |f(y)| \, dy
\end{equation}
for every $\lambda > 0$, where $|E_\lambda|$ is the Lebesgue measure
of $E_\lambda$.  In particular, $f^* < \infty$ almost everywhere on
${\bf R}$.  To prove this estimate, it suffices to show that
\begin{equation}
        |K| \le \frac{2}{\lambda} \int_{\bf R} |f(y)| \, dy
\end{equation}
for every compact set $K \subseteq E_\lambda$.

        By compactness, there are finitely many open intervals
$I_1, \ldots, I_n$ such that
\begin{equation}
        \frac{1}{|I_l|} \int_{I_l} |f(y)| \, dy > \lambda
\end{equation}
for each $l$, and
\begin{equation}
        K \subseteq \bigcup_{l = 1}^n I_l.
\end{equation}
We may also suppose that each $x \in {\bf R}$ is contained in at most
two of the $I_l$'s.  For if a point is contained in three of these
intervals, then one of the intervals is contained in the union of the
other two, and may be dropped from the collection.  This implies that
\begin{eqnarray}
|K| \le \sum_l |I_l| & \le & \sum_l \frac{1}{\lambda} \int_{I_l} |f(y)| \, dy\\
 & \le & \frac{2}{\lambda} \int_{\bigcup_l I_l} |f(y)| \, dy \nonumber\\
 & \le & \frac{2}{\lambda} \int_{\bf R} |f(y)| \, dy, \nonumber
\end{eqnarray}
as desired.

\section{$L^p$ estimates}
\label{L^p estimates}
\setcounter{equation}{0}

        Let $\phi$ be a nonnegative measurable function on the real
line, and let $p$ be a real number with $p \ge 1$.  Observe that
\begin{eqnarray}
 \int_{\bf R} \phi(x)^p \, dx & = &
 \int_{\bf R} \int_0^{\phi(x)} p \, \lambda^{p - 1} \, d\lambda \, dx \\
  & = & \int_0^\infty p \, \lambda^{p - 1} \,
         |\{x \in {\bf R} : \phi(x) > \lambda\}| \, d\lambda. \nonumber
\end{eqnarray}
We would like to apply this to $\phi = f^*$, where $f \in L^p({\bf
R})$ and $p > 1$.

        For each $\lambda > 0$, let $f_\lambda$ be the function
defined by
\begin{eqnarray}
 f_\lambda(x) & = & f(x) \quad\hbox{when } |f(x)| \le \frac{\lambda}{2} \\
              & = & 0    \enspace\qquad\hbox{otherwise}. \nonumber
\end{eqnarray}
Thus
\begin{equation}
        f^*(x) \le f_\lambda^*(x) + (f - f_\lambda)^*(x)
             \le \frac{\lambda}{2} + (f - f_\lambda)(x)^*
\end{equation}
for every $x \in {\bf R}$, by subadditivity of $f^*$ and
(\ref{||f^*||_infty le ||f||_infty}), which implies that
\begin{equation}
 \{x \in {\bf R} : f^*(x) > \lambda\} \subseteq
  \bigg\{x \in {\bf R} : (f - f_\lambda)^*(x) > \frac{\lambda}{2}\bigg\}.
\end{equation}

        Combining this with (\ref{|E_lambda| le frac{2}{lambda}
int_{bf R} |f(y)| dy}) applied to $f - f_\lambda$, we get that
\begin{equation}
        |\{x \in {\bf R} : f^*(x) > \lambda\}| \le
 \frac{4}{\lambda} \int_{\{y \in {\bf R} : |f(y)| > \lambda/2\}} |f(y)| \, dy
\end{equation}
for every $\lambda > 0$.  Hence
\begin{equation}
 \int_{\bf R} f^*(x)^p \, dx \le \int_0^\infty 4 \, p \, \lambda^{p - 2} \,
       \int_{\{y \in {\bf R} : |f(y)| > \lambda/2\}} |f(y)| \, dy \, d\lambda.
\end{equation}
Interchanging the order of integration, we get that
\begin{equation}
 \int_{\bf R} f^*(x)^p \, dx \le \int_{\bf R} \int_0^{2 \, |f(y)|} 4 \, p \,
                              \lambda^{p - 2} \, |f(y)| \, d\lambda \, dy.
\end{equation}
If $p > 1$, then it follows that
\begin{equation}
        \int_{\bf R} f^*(x)^p \, dx
         \le \frac{4 \, p \, 2^{p - 1}}{p - 1} \int_{\bf R} |f(y)|^p \, dy.
\end{equation}

        If $f$ is not equal to $0$ almost everywhere on ${\bf R}$,
then one can check that $f^*(x)$ is at least a positive constant
multiple of $1/|x|$ when $|x|$ is sufficiently large.  This implies
that $f^*$ is not integrable on ${\bf R}$, although it can still be
locally integrable.

\section{Maximal functions, 2}
\label{maximal functions, 2}
\setcounter{equation}{0}

        Let $(M, d(x, y))$ be a metric space.  For each $x \in M$ and
$r > 0$, the open ball with center $x$ and radius $r$ is defined by
\begin{equation}
        B(x, r) = \{y \in M : d(x, y) < r \}.
\end{equation}
If $B = B(x, r)$ is an open ball and $a > 0$, then $a \, B$ denotes
the ball $B(x, a \, r)$.

        Suppose that $B_1, B_2, \ldots$ is a finite sequence of open
balls in $M$, or an infinite sequence of balls whose radii converge to
$0$.  Let $\beta_1$ be one of these balls with maximum radius, and let
$\beta_2$ be another one of these balls which is disjoint from
$\beta_1$ and has maximal radius, if there is one.  Continuing in this
manner, we get a finite or infinite sequence of pairwise disjoint
balls $\beta_1, \beta_2, \ldots$ chosen from $B_1, B_2, \ldots$ such
that for each $B_i$ there is a $\beta_j$ which intersects $B_i$ and
has radius at least that of $B_i$.  Hence $B_i \subseteq 3 \,
\beta_j$, and so
\begin{equation}
        \bigcup_i B_i \subseteq \bigcup_j 3 \beta_j.
\end{equation}

        A positive Borel measure $\mu$ on $M$ is said to be
\emph{doubling} if the measure of each open ball in $M$ is positive
and finite, and if there is a $C > 0$ such that
\begin{equation}
        \mu(2 \, B) \le C \, \mu(B)
\end{equation}
for every open ball $B$ in $M$.  If $f$ is a locally integrable
function on $M$ with respect to $\mu$, then the corresponding maximal
function $f^*$ is defined by
\begin{equation}
        f^*(x) = \sup_{x \in B} \frac{1}{\mu(B)} \int_B |f| \, d\mu,
\end{equation}
where the supremum is taken over all open balls $B$ in $M$ with $x \in B$.
As in \cite{c-w-1, c-w-2}, one can show that $f^*$ satisfies properties
like those in the previous sections, using the covering lemma described in
the preceding paragraph in particular.

        For example, this applies to Lebesgue measure on ${\bf R}^n$
equipped with the standard metric.  As in \cite{j1}, one can also
consider ${\bf R}^{n + 1}$ equipped with a metric associated to
parabolic dilations.  In addition to nonstandard dilations, there can
be more complicated geometric structure related to nilpotent Lie
groups.  Many common fractals are equipped with natural doubling
measures.

\section{Discrete maximal functions}
\label{discrete maximal functions}
\setcounter{equation}{0}

        If $f$ is a function on the set ${\bf Z}$ of integers,
then the corresponding maximal function is defined by
\begin{equation}
 f^*(l) = \sup_{a \le l \le b} \frac{1}{b - a + 1} \sum_{j = a}^b |f(j)|,
\end{equation}
where more precisely the supremum is taken over all integers $a$, $b$
with $a \le l \le b$, and may be infinite.  As before, $f \mapsto f^*$
is sublinear, $f^*$ is bounded when $f$ is bounded, and satisfies
\begin{equation}
        \sup_{l \in {\bf Z}} f^*(l) \le \sup_{l \in {\bf Z}} |f(l)|.
\end{equation}
If we identify functions on ${\bf Z}$ with functions on ${\bf R}$ that
are constant on intervals of the form $[l, l + 1)$ for $l \in {\bf Z}$
and have the same values at the integers, then this maximal function
is obviously less than or equal to the previous one on ${\bf R}$.
Using this observation, or arguments analogous to the earlier ones,
it is easy to see that
\begin{equation}
        |\{j \in {\bf Z} : f^*(j) > \lambda\}|
         \le \frac{2}{\lambda} \sum_{j = -\infty}^\infty |f(j)|
\end{equation}
for every $\lambda > 0$, where $|E|$ is the number of elements of $E
\subseteq {\bf Z}$, and
\begin{equation}
        \sum_{j = -\infty}^\infty f^*(j)^p
 \le \frac{4 \, p \, 2^{p - 1}}{p - 1} \sum_{j = -\infty}^\infty |f(j)|^p
\end{equation}
when $1 < p < \infty$.

\section{Another variant}
\label{another variant}
\setcounter{equation}{0}

        Let $X$ be a set with a positive measure $\mu$, and consider
$X \times {\bf Z}$ with the measure which is the product of $\mu$ and
counting measure.  If $f$ is a measurable function on $X \times {\bf
Z}$, which is to say that $f(x, l)$ is measurable in $x$ for each $l$,
then consider the maximal function $f^*$ defined by
\begin{equation}
 f^*(x, l) = \sup_{a \le l \le b} \frac{1}{b - a + 1} \sum_{j = a}^b |f(x, j)|,
\end{equation}
which is the maximal function of $f(x, l)$ in $l$ for each $x \in X$.
Thus $f \mapsto f^*$ is sublinear, $f^*$ is bounded when $f$ is
bounded, and
\begin{equation}
        \|f\|_\infty \le \|f\|_\infty.
\end{equation}
Moreover,
\begin{equation}
 \sum_{j = -\infty}^\infty \mu(\{x \in X : f^*(x, j) > \lambda\})
 \le \frac{2}{\lambda} \sum_{j = -\infty}^\infty \int_X |f(x, j)| \, d\mu(x)
\end{equation}
for each $\lambda > 0$, and
\begin{equation}
 \sum_{j = -\infty}^\infty \int_X f^*(x, j)^p \, d\mu(x)
 \le \frac{4 \, p \, 2^{p - 1}}{p - 1}
      \sum_{j = -\infty}^\infty \int_X |f(x, j)|^p \, d\mu(x)
\end{equation}
when $1 < p < \infty$.

\section{Transference}
\label{transference}
\setcounter{equation}{0}

        Let $X$ be a set with a positive measure $\mu$, let $\phi$ be
a one-to-one measure-preserving mapping of $X$ onto itself, and let
\begin{equation}
        T(f) = f \circ \phi
\end{equation}
be the corresponding linear transformation on measurable functions on
$X$.  Consider the maximal function
\begin{equation}
 A^*(f)(x) = \sup_{k, l \ge 0} \frac{1}{k + l + 1} \sum_{j = -k}^l |T^j(f)(x)|,
\end{equation}
where the supremum is taken over nonnegative integers $k$, $l$.
It is also convenient to consider the approximations to this defined by
\begin{equation}
 A_n^*(f)(x) = \max_{0 \le k, l \le n} \frac{1}{k + l + 1}
                                        \sum_{j = -k}^l |T_j(f)(x)|.
\end{equation}
Thus $A_n^*(f)(x)$ is monotone increasing in $n$, and tends to
$A^*(f)(x)$ as $n \to \infty$.  As usual, $A^*$ and $A_n^*$ are
sublinear operators that send bounded functions to bounded functions, with
\begin{equation}
        \|A_n^*(f)\|_\infty \le \|A^*(f)\|_\infty \le \|f\|_\infty.
\end{equation}
By construction, they also commute with $T$, which is to say that
\begin{equation}
        A_n^*(T(f)) = T(A_n^*(f)), \quad A^*(T(f)) = T(A^*(f)).
\end{equation}
In particular,
\begin{equation}
        \mu(\{x \in X : A_n^*(T^l(f))(x) > \lambda\})
         = \mu(\{x \in X : A_n^*(f)(x) > \lambda\})
\end{equation}
for each $n$, $l$, and $\lambda$, and
\begin{equation}
        \int_X A_n^*(T^l(f))^p \, d\mu = \int_X A_n^*(f)^p \, d\mu
\end{equation}
for each $n$, $l$, and $p$.

        Let $N$ be a large positive integer.  If $f$ is a measurable
function on $X$, then let $F(x, j)$ be the function on $X \times {\bf
Z}$ defined by
\begin{equation}
        F(x, j) = T^j(f)(x)
\end{equation}
when $0 \le j \le N$, and $F(x, j) = 0$ otherwise.  Because $\phi$ is
measure-preserving,
\begin{equation}
        \sum_{j = -\infty}^\infty \int_X |F(x, j)|^p \, d\mu(x)
                        = (N + 1) \int_X |f|^p \, d\mu
\end{equation}
for each $p$.  Let $F^*(x, l)$ be the maximal function in the second
variable, as in the previous section.  If $n \le l \le N - n$, then
\begin{equation}
        A_n^*(T^l(f))(x) \le F^*(x, l).
\end{equation}
It follows that
\begin{eqnarray}
\lefteqn{(N - 2 \, n) \, \mu(\{x \in X : A_n^*(f)(x) > \lambda\})} \\
 & & \le \sum_{j = -\infty}^\infty \mu(\{x \in X : F^*(x, j) > \lambda\})
                                                                  \nonumber
\end{eqnarray}
for each $\lambda$, and hence
\begin{eqnarray}
\lefteqn{(N - 2 \, n) \, \mu(\{x \in X : A_n^*(f)(x) > \lambda\})} \\
& & \le \frac{2}{\lambda} \sum_{j = -\infty}^\infty \int_X |F(x, j)| \, d\mu(x)
                                                    \nonumber \\
  & & = (N + 1) \frac{2}{\lambda} \int_X |f| \, d\mu. \nonumber
\end{eqnarray}
Equivalently,
\begin{equation}
        \mu(\{x \in X : A_n^*(f)(x) > \lambda\})
 \le \frac{(N + 1)}{(N - 2 \, n)} \, \frac{2}{\lambda} \int_X |f| \, d\mu
\end{equation}
for every $N > 2 \, n$, and so
\begin{equation}
        \mu(\{x \in X : A_n^*(f)(x) > \lambda\})
         \le \frac{2}{\lambda} \int_X |f| \, d\mu,
\end{equation}
which implies that
\begin{equation}
 \mu(\{x \in X : A^*(f)(x) > \lambda\}) \le \frac{2}{\lambda} \int_X |f| d\,\mu
\end{equation}
for every $\lambda$.

        Similarly, for each $p$,
\begin{equation}
        (N - 2 \, n) \int_X A_n^*(f)^p \, d\mu
         \le \sum_{j = -\infty}^\infty \int_X F^*(x, j)^p \, d\mu(x).
\end{equation}
If $1 < p < \infty$, then
\begin{eqnarray}
\lefteqn{(N - 2 \, n) \int_X A_n^*(f)^p \, d\mu} \\
 & & \le \frac{4 \, p \, 2^{p - 1}}{p - 1}
  \sum_{j = -\infty}^\infty \int_X |F(x, j)|^p \, d\mu(x) \nonumber \\
& & = (N + 1) \, \frac{4 \, p \, 2^{p-1}}{p - 1} \int_X |f|^p \, d\mu.\nonumber
\end{eqnarray}
Thus
\begin{equation}
        \int_X A_n^*(f)^p \, d\mu
          \le \frac{(N + 1)}{(N - 2 \, n)} \,
           \frac{4 \, p \, 2^{p - 1}}{p - 1} \int_X |f|^p \, d\mu
\end{equation}
when $N > 2 \, n$, which implies that
\begin{equation}
        \int_X A_n^*(f)^p \, d\mu
         \le \frac{4 \, p \, 2^{p - 1}}{p - 1} \int_X |f|^p \, d\mu,
\end{equation}
and hence
\begin{equation}
        \int_X A^*(f)^p \, d\mu
         \le \frac{4 \, p \, 2^{p - 1}}{p - 1} \int_X |f|^p \, d\mu.
\end{equation}

\section{Almost-everywhere convergence}
\label{almost-everywhere convergence}
\setcounter{equation}{0}

        A standard application of maximal function estimates is to the
existence of pointwise limits almost everywhere.  For example, if $f$
is a locally integrable function on the real line, then
\begin{equation}
        \frac{1}{|I|} \int_I |f(y) - f(x)| \, dy \to 0
\end{equation}
for almost every $x \in {\bf R}$, where the limit is taken for open
intervals $I$ with $x \in I$ as $|I| \to 0$.  One may as well suppose
that $f$ is integrable here, or even has compact support, since the
problem is local.  The limit is trivial when $f$ is continuous, and
any integrable function can be approximated in the $L^1$ norm by a
continuous function.  Using maximal function estimates, one can show
that the approximations by continuous functions also work for the
pointwise limits almost everywhere.

        One can also show that maximal functions associated to the
Abel and Cesaro sums of the Fourier series of an integrable function
are bounded by the analogue of the Hardy--Littlewood maximal function
on the circle.  Hence these maximal functions satisfy the same sort of
estimates as before.  This can be used to show that the Abel and
Cesaro sums converge to the function almost everywhere, which could be
derived from the limit described in the previous paragraph as well,
and which is basically the same type of statement anyway.

        Suppose now that $X$ is a set with a positive measure $\mu$,
$\phi$ is a measure-preserving transformation on $X$, and that
$T(f) = f \circ \phi$ for measurable functions $f$ on $X$.  
If $f \in L^p(X)$, $1 < p < \infty$, then
\begin{eqnarray}
 \int_X \sum_{n = 0}^\infty \frac{|T^n(f)(x)|^p}{(n + 1)^p} \, d\mu(x)
 & = & \sum_{n = 0}^\infty \frac{1}{(n + 1)^p} \int_X |T^n(f)(x)|^p \, d\mu \\
 & = & \sum_{n = 0}^\infty \frac{1}{(n + 1)^p} \int_X |f|^p \, d\mu < \infty.
                                                                   \nonumber
\end{eqnarray}
Hence
\begin{equation}
        \sum_{n = 0}^\infty \frac{|T^n(f)(x)|^p}{(n + 1)^p} < \infty
\end{equation}
for almost every $x \in X$, and thus
\begin{equation}
        \lim_{n \to \infty} \frac{T^n(f)(x)}{n + 1} = 0
\end{equation}
for almost every $x \in X$.  Of course,
\begin{equation}
        \frac{\|T^n(f)\|_\infty}{n + 1} = \frac{\|f\|_\infty}{n + 1} \to 0
\end{equation}
as $n \to \infty$ when $f$ is bounded $X$.

        For each $x \in X$, consider the sequence of averages
\begin{equation}
\label{1 / n + 1 sum_{j = 0}^n T^j(f)(x)}
        \frac{1}{n + 1} \sum_{j = 0}^n T^j(f)(x).
\end{equation}
If $T(f) = f$, then
\begin{equation}
        \frac{1}{n + 1} \sum_{j = 0}^n T^j(f)(x) = f(x)
\end{equation}
for each $n$, and converge as $n \to \infty$ trivially.  If $f = T(b)
- b$ for some $b \in L^p(X)$, $1 < p \le \infty$, then
\begin{equation}
\frac{1}{n + 1} \sum_{j = 0}^n T^j(f)(x) = \frac{T^{n + 1}(b)(x) - b(x)}{n + 1}
\end{equation}
converges to $0$ almost everywhere on $X$ as $n \to \infty$, by the
remarks in the preceding paragraph.  As in Section \ref{unitary
transformations}, the functions of the form
\begin{equation}
        a + (T(b) - b), \ a, b \in L^2(X), \, T(a) = a,
\end{equation}
are dense in $L^2(X)$, since $T$ defines a unitary transformation on
$L^2(X)$.  Thus the averages (\ref{1 / n + 1 sum_{j = 0}^n T^j(f)(x)})
converge almost everywhere on $X$ for a dense class of functions in
$L^2(X)$.  Using the boundedness of the maximal function on $L^2(X)$,
one can show that the averages converge almost everywhere for any $f
\in L^2(X)$.  This also holds for $f \in L^p(X)$ when $1 \le p <
\infty$, because $L^p(X) \cap L^2(X)$ is dense in $L^p(X)$, and using
the maximal function estimates on $L^p(X)$.

\section{Multiplication operators}
\label{multiplication operators}
\setcounter{equation}{0}

        Let $(X, \mu)$ be a measure space, and let $b$ be a bounded
measurable function on $X$.  The corresponding multiplication operator
is defined by
\begin{equation}
        T(f) = b \, f.
\end{equation}
This determines a bounded linear operator on $L^p(X)$ for $1 \le p \le
\infty$, with operator norm equal to the essential supremum norm of
$b$ for each $p$.

        In particular, if $\{b_j\}_j$ is a sequence of bounded
measurable functions on $X$ that converges to $b$ in the $L^\infty$
norm, then the corresponding sequence $\{T_j\}_j$ of multiplication
operators converges to $T$ in the operator norm on $L^p(X)$ for each
$p$, $1 \le p \le \infty$.  If instead $\{b_j\}_j$ is a sequence of
bounded measurable functions with uniformly bounded $L^\infty$ norms
that converges to $b$ pointwise almost everywhere on $X$, then the
dominated convergence theorem implies that $\{T_j(f)\}_j$ converges to
$T(f)$ in the $L^p$ norm for every $f \in L^p(X)$ when $1 \le p < \infty$.

        For example, if $\|b\|_\infty < 1$, then $b^j$ converges to
$0$ in the $L^\infty$ norm as $j \to \infty$.  If $\|b\|_\infty = 1$
and $|b(x)| < 1$ for almost every $x \in X$, then $\|b^j\|_\infty = 1$
for each $j$ and $\lim_{j \to \infty} b(x)^j = 0$ almost everywhere on
$X$.

        Suppose that $\|b\|_\infty \le 1$, and consider
\begin{equation}
        a_n = \frac{1}{n + 1} \sum_{j = 0}^n b^j.
\end{equation}
Thus
\begin{equation}
        \frac{1}{n + 1} \sum_{j = 0}^n T^j(f) = a_n \, f.
\end{equation}
Clearly $\|a_n\|_\infty \le 1$ for each $n$.  As in Section
\ref{cesaro means}, $\{a_n\}_n$ converges pointwise to the function
that is equal to $1$ when $b$ is equal to $1$ and to $0$ elsewhere.
Hence the sequence of averages of the $T_j$'s converges strongly to
the operator of multiplication by the characteristic function of the
set where $b = 1$ on $L^p(X)$ when $1 \le p < \infty$.

\section{Doubling spaces}
\label{doubling spaces}
\setcounter{equation}{0}

        A metric space $(M, d(x, y))$ is said to be \emph{doubling} if
there is a $C > 0$ such that every ball in $M$ of radius $r$ is
contained in the union of $\le C$ balls of radius $r/2$.  Thus every
ball of radius $r$ can also be covered by $\le C^2$ balls of radius
$r/4$, etc.  Equivalently, $M$ is doubling if there is a $C' > 0$ such
that every set in $M$ of diameter $t$ can be covered by $\le C'$
subsets of diameter $\le t/2$.

        For example, it is easy to see that finite-dimensional
Euclidean spaces are doubling.  In this case, it suffices to consider
coverings of the unit ball, because a covering of any other ball can
be reduced to the unit ball using translations and dilations.  A
similar argument works on nilpotent Lie groups under suitable
conditions.  If $(M, d(x, y))$ is a metric space which is doubling and
$E \subseteq M$, then $E$ is doubling with respect to the restriction
of $d(x, y)$ to $E$.  In particular, subsets of finite-dimensional
Euclidean spaces are doubling with respect to the induced metric.

        Let $M$ be a metric space with a doubling measure $\mu$, and
let $B$ be a ball of radius $r$ in $M$.  Suppose that $A \subseteq B$
has the property that $d(x, y) \ge r/2$ for every $x, y \in A$ with $x
\ne y$.  Thus the balls $B(x, r/2)$, $x \in A$, are pairwise disjoint,
and $B(x, r/2) \subseteq (3/2) B$ for every $x \in B$, and hence
\begin{equation}
        \sum_{x \in A} \mu(B(x, r/2)) \le \mu((3/2)B).
\end{equation}
This implies an upper bound for the number of elements of $A$ in terms
of the doubling constant for $\mu$, since
\begin{equation}
        B \subseteq B(x, 2 \, r)
\end{equation}
for each $x \in B$, so that $\mu(B)$ is less than or equal to a
constant multiple of $\mu(B(x, r/2))$.  If $A$ is a maximal set of
this type, then for each $z \in B$ there is an $x \in A$ such that
$d(x, z) \le r/2$, which means that
\begin{equation}
        B \subseteq \bigcup_{x \in A} B(x, r/2),
\end{equation}
from which it follows that $M$ is doubling.

        Remember that a set $E$ is a metric space $M$ is said to be
\emph{totally bounded} if for each $\epsilon > 0$, $E$ can be covered
by finitely many balls of radius $\epsilon$.  Totally bounded subsets
of any metric space are automatically bounded, and bounded subsets of
a doubling metric space are totally bounded.  If $M$ is complete, then
$K \subseteq M$ is compact if and only if $K$ is closed and totally
bounded.  If $M$ is also doubling, then it follows that $K \subseteq
M$ is compact if and only if $K$ is closed and bounded.

\section{Ultrametrics}
\label{ultrametrics}
\setcounter{equation}{0}

        A metric $d(x, y)$ on a set $M$ is said to be an
\emph{ultrametric} if
\begin{equation}
        d(x, z) \le \max (d(x, y), d(y, z))
\end{equation}
for every $x, y, z \in M$, which is stronger than the usual triangle
inequality.  Some examples of ultrametric spaces are given in the next section.

        Let $(M, d(x, y))$ be an ultrametric space.  If $x, y \in M$
and $r, t > 0$, then either
\begin{equation}
        B(x, r) \subseteq B(y, t),
\end{equation}
or
\begin{equation}
        B(y, t) \subseteq B(x, r),
\end{equation}
or
\begin{equation}
        B(x, r) \cap B(y, t) = \emptyset.
\end{equation}
More precisely, the first alternative holds when $d(x, y)$ and $r$ are
less than or equal to $t$, the second alternative holds when $d(x, y)$
and $t$ are less than or equal to $r$, and the third alternative holds
when $r$ and $t$ are strictly less than $d(x, y)$.

        Suppose that $\mathcal{B}$ is a collection of open balls in
$M$.  The remarks of the previous paragraph imply that the maximal
elements of $\mathcal{B}$ are pairwise disjoint.  Under suitable
conditions, every element of $\mathcal{B}$ is contained in a maximal
element.

        This is a very simple covering argument for ultrametric
spaces.  Doubling conditions are still important, even if they are not
required for a covering lemma.

\section{Sequence spaces, 2}
\label{sequence spaces, 2}
\setcounter{equation}{0}

        Let $A$ be a set with at least two elements, and let $M$ be
the set of sequences $x = \{x_j\}_{j = 1}^\infty$ with $x_j \in A$ for
each $j$.  Fix a positive real number $\rho < 1$, and put
\begin{equation}
        d(x, y) = \rho^n
\end{equation}
when $x, y \in M$, $x \ne y$, and $n$ is the largest nonnegative
integer such that $x_j = y_j$ for $j \le n$, and otherwise put $d(x,
y) = 0$.  This defines an ultrametric on $M$ for which the
corresponding topology is the product topology, when $M$ is considered
as the product of infinitely many copies of $A$ equipped with the
discrete topology.  Let us restrict our attention to the case where
$A$ has only finitely many elements, so that $M$ is compact.  It is
easy to see that $M$ is doubling in this case.  Let $w(a)$ be a
positive real number for each $a \in A$ whose sum is equal to $1$.
This defines a probability measure on $A$ which leads to a probability
measure $\mu$ on $M$, and one can check that $\mu$ is a doubling
measure on $M$.

\section{Snowflake transforms}
\label{snowflakes}
\setcounter{equation}{0}

        If $a$, $x$, $y$ are nonnegative real numbers and $0 < a < 1$, then
\begin{equation}
        (x + y)^a \le x^a + y^a.
\end{equation}
This is because
\begin{equation}
        x + y \le (x^a + y^a) \, \max(x, y)^{1 - a}
\end{equation}
and
\begin{equation}
        \max(x, y) \le (x^a + y^a)^{1/a}
\end{equation}
imply that $x + y \le (x^a + y^a)^{1/a}$.

        If $(M, d(x, y))$ is a metric space and $0 < a < 1$, then it
follows that $d(x, y)^a$ is also a metric on $M$.  If $d(x, y)$ is an
ultrametric on $M$, then $d(x, y)^a$ is an ultrametric on $M$ for
every $a > 0$.  Note that $d(x, y)^a$ determines the same topology on
$M$ as $d(x, y)$ in both situations.

        For each $p \in M$ and $r > 0$, the open ball in $M$ centered
at $p$ with radius $r$ with respect to $d(x, y)$ is the same as the
open ball centered at $p$ with radius $r^a$ with respect to $d(x,
y)^a$.  Using this fact, it is easy to see that $M$ is doubling with
respect to $d(x, y)$ if and only if it is doubling with respect to
$d(x, y)^a$.  Similarly, a positive Borel measure $\mu$ on $X$ is
doubling with respect to $d(x, y)$ if and only if it is doubling with
respect to $d(x, y)^a$.  Of course, the doubling constants for $d(x,
y)^a$ are not normally the same as for $d(x, y)$.

        The diameter of a set $E \subseteq M$ with respect to $d(x,
y)^a$ is equal to the $a$th power of the diameter of $E$ with respect
to $d(x, y)$.  It is easy to see that the $t$-dimensional Hausdorff
measure of $E \subseteq M$ with respect to $d(x, y)$ is equal to the
$(t/a)$-dimensional Hausdorff measure of $E$ with respect to $d(x,
y)^a$, using the preceding observation and the definition of Hausdorff
measures.  In particular, the Hausdorff dimension of $E$ with respect
to $d(x, y)^a$ is equal to the Hausdorff dimension of $E$ with respect
to $d(x, y)$ divided by $a$.

\section{Invariant measures}
\label{invariant measures}
\setcounter{equation}{0}

        Let $X$ be a compact Hausdorff topological space, and let
$\mathcal{C}(X)$ be the space of continuous real-valued functions on
$X$ equipped with the supremum norm.  A linear functional $\lambda$ on
$\mathcal{C}(X)$ is said to be \emph{bounded} if it is bounded as a
linear mapping into ${\bf R}$ as a one-dimensional vector space with
the absolute value as norm.  Thus $\lambda$ is bounded if there is an
$A \ge 0$ such that
\begin{equation}
        |\lambda(f)| \le A \, \|f\|_{sup}
\end{equation}
for every $f \in \mathcal{C}(X)$, in which case the dual norm $\|\lambda\|_*$
of $\lambda$ is the smallest value of $A$ for which this condition holds.
This defines a norm on the vector space $\mathcal{C}(X)^*$ of bounded
linear functionals on $\mathcal{C}(X)$, which is the same as the operator
norm of $\lambda$ as a linear mapping from $\mathcal{C}(X)$ into ${\bf R}$.

        A linear functional $\lambda$ on $\mathcal{C}(X)$ is said to be
\emph{nonnegative} if
\begin{equation}
        \lambda(f) \ge 0
\end{equation}
for every $f \in \mathcal{C}(X)$ which is nonnegative in the sense
that $f(x) \ge 0$ for every $x \in X$.  If $\lambda$ is a nonnegative
linear functional on $\mathcal{C}(X)$, then $\lambda$ is bounded, and
\begin{equation}
        \|\lambda\|_* = \lambda(1),
\end{equation}
where the right side means $\lambda$ applied to the constant function
equal $1$ on $X$.  Every nonnegative linear functional on
$\mathcal{C}(X)$ can be expressed by integration with respect to a
positive Borel measure on $X$, by the Riesz representation theorem.

        A set $U \subseteq \mathcal{C}(X)^*$ is an open set in the
weak$^*$ topology if for every $\lambda \in U$ there are finitely
many continuous functions $f_1, \ldots, f_n$ on $X$ and positive real
numbers $r_1, \ldots r_n$ such that
\begin{equation}
 \ \{\lambda' \in \mathcal{C}(X)^* : |\lambda'(f_1) - \lambda(f_1)| < r_1,
                            \ldots, |\lambda'(f_n) - \lambda(f_n)| < r_n \}
                                                           \subseteq U.
\end{equation}
The Banach--Alaoglu theorem implies that
\begin{equation}
        B^* = \{\lambda \in \mathcal{C}(X)^* : \|\lambda\|_* \le 1\}
\end{equation}
is compact with respect to this topology.  If the topology on $X$ is
determined by a metric, then $\mathcal{C}(X)$ is separable, and one
can show that the topology on $B^*$ induced by the weak$^*$ topology
is determined by a metric.

        Let $\phi$ be a homeomorphism of $X$ onto itself, and let
\begin{equation}
        T(f) = f \circ \phi
\end{equation}
be the corresponding linear operator on $\mathcal{C}(X)$.  This leads
to a dual linear operator $T^*$ on $\mathcal{C}(X)^*$, defined by
\begin{equation}
        T^*(\lambda)(f) = \lambda(T(f)).
\end{equation}
Note that $T(f) \ge 0$ when $f \ge 0$, and hence $T^*(\lambda)$ is
nonnegative when $\lambda$ is nonnegative.  Similarly, $T$ preserves
the supremum norm on $\mathcal{C}(X)$, and so $T^*$ preserves the dual
norm on $\mathcal{C}(X)^*$.  Of course, $T$ automatically sends
constant functions to themselves on $X$.

        Consider
\begin{equation}
        \lambda_n = \frac{1}{n + 1} \sum_{j = 0}^n (T^*)^j(\lambda)
\end{equation}
for each nonnegative integer $n$ and $\lambda \in \mathcal{C}(X)^*$.
Observe that
\begin{equation}
        T^*(\lambda_n) - \lambda_n =
           \frac{1}{n + 1}((T^*)^{n + 1} (\lambda) - \lambda),
\end{equation}
which implies that
\begin{equation}
        \|T^*(\lambda_n) - \lambda_n\|_* \le \frac{2 \, \|\lambda\|_*}{n + 1}.
\end{equation}
In particular,
\begin{equation}
 \lim_{n \to \infty} \|T^*(\lambda_n) - \lambda_n\|_* = 0
\end{equation}

        Let $E_l$ be the closure of $\lambda_l, \lambda_{l + 1},
\ldots$ in the weak$^*$ topology.  By compactness,
\begin{equation}
        E = \bigcap_{l = 0}^\infty E_l \ne \emptyset.
\end{equation}
If the topology on $X$ is determined by a metric, then there is a
subsequence of $\lambda_0, \lambda_1, \ldots$ that converges in the
weak$^*$ topology, and whose limit is in $E$.  Every element of $E$ is
invariant under $T^*$, by the remarks of the previous paragraph.

        If $\lambda$ is nonnegative, then $\lambda_n$ is nonnegative
for each $n$, and every element of $E$ is nonnegative.  In this case,
$\lambda_n(1) = \lambda(1)$ for each $n$, and every element of $E$ has
the same value when applied to the constant function $1$.  It follows
that $\lambda_n$ and every element of $E$ has the same norm as
$\lambda$ when $\lambda \ge 0$.  Using this, one can check that there
is a nonnegative linear functional on $\mathcal{C}(X)$ with norm $1$
which is invariant under $T^*$, and which corresponds to an invariant
probability measure on $X$.


\begin{thebibliography}{99}

\addcontentsline{toc}{section}{References}

\bibitem {a} J.~Aaronson, {\it An Introduction to Infinite Ergodic
Theory}, American Mathematical Society, 1997.

\bibitem {a-t} L.~Ambrosio and P.~Tilli, {\it Topics on Analysis in
Metric Spaces}, Oxford University Press, 2004.

\bibitem {arv} W.~Arveson, {\it A Short Course on Spectral Theory},
Springer-Verlag, 2002.

\bibitem {b1} R.~Beals, {\it Topics in Operator Theory}, University of
Chicago Press, 1971.

\bibitem {b2} R.~Beals, {\it Analysis: An Introduction}, Cambridge
University Press, 2004.

\bibitem {b-c} S.~Bochner and K.~Chandrasekharan, {\it Fourier
Transforms}, Princeton University Press, 1949.

\bibitem {c-w-1} R.~Coifman and G.~Weiss, {\it Analyse Harmonique
Non-Commutative sur certains Espaces Homog\`enes}, Lecture Notes in
Mathematics {\bf 242}, Springer-Verlag, 1971.

\bibitem {c-w-2} R.~Coifman and G.~Weiss, {\it Transference Methods in
Analysis}, American Mathematical Society, 1976.

\bibitem {c-w-3} R.~Coifman and G.~Weiss, {\it Extensions of Hardy
spaces and their use in analysis}, Bulletin of the American
Mathematical Society {\bf 83} (1977), 569--645.

\bibitem {cw} J.~Conway, {\it A Course in Functional Analysis}, 2nd
edition, Springer-Verlag, 1990.

\bibitem {dav-s} G.~David and S.~Semmes, {\it Fracture Fractals and
Broken Dreams: Self-Similar Geometry through Metric and Measure},
Oxford University Press, 1997.

\bibitem {dev-s} R.~DeVore and R.~Sharpley, {\it Maximal Functions
Measuring Smoothness}, Memoirs of the American Mathematical Society
{\bf 293}, 1984.

\bibitem {duo} J.~Duoandikoetxea, {\it Fourier Analysis}, translated
and revised from the 1985 Spanish original by D.~Cruz-Uribe, American
Mathematical Society, 2001.

\bibitem {e-m-t} Y.~Eidelman, V.~Milman, and A.~Tsolomitis, {\it
Functional Analysis: An Introduction}, American Mathemtical Society,
2004.

\bibitem {fa1} K.~Falconer, {\it The Geometry of Fractal Sets},
Cambridge University Press, 1986.

\bibitem {fa2} K.~Falconer, {\it Techniques in Fractal Geometry},
Wiley, 1997.

\bibitem {fa3} K.~Falconer, {\it Fractal Geometry: Mathematical
Foundations and Applications}, 2nd edition, Wiley, 2003.

\bibitem {fe} H.~Federer, {\it Geometric Measure Theory},
Springer-Verlag, 1969.

\bibitem {fo} G.~Folland, {\it Real Analysis}, 2nd edition, Wiley,
1999.

\bibitem {f-s} G.~Folland and E.~Stein, {\it Hardy Spaces on
Homogeneous Groups}, Princeton University Press, 1982.

\bibitem {hf} H.~Furstenberg, {\it Recurrence in Ergodic Theory and
Combinatorial Number Theory}, Princeton University Press, 1981.

\bibitem {gc-rf} J.~Garc\'{\i}a-Cuerva and J.~Rubio de Francia, {\it
Weighted Norm Inequalities and Related Topics}, North-Holland, 1984.

\bibitem {g} R.~Goldberg, {\it Methods of Real Analysis}, 2nd edition,
Wiley, 1976.

\bibitem {h1} P.~Halmos, {\it Lectures on Ergodic Theory}, Chelsea,
1960.

\bibitem {h2} P.~Halmos, {\it Finite-Dimensional Vector Spaces},
Springer-Verlag, 1974.

\bibitem {h3} P.~Halmos, {\it A Hilbert Space Problem Book}, 2nd
edition, Springer-Verlag, 1982.

\bibitem {h4} P.~Halmos, {\it Introduction to Hilbert Space and the
Theory of Spectral Multiplicity}, 2nd edition, AMS Chelsea Publishing,
1998.

\bibitem {h} V.~Hansen, {\it Functional Analysis: Entering Hilbert
Space}, World Scientific, 2006.

\bibitem {jh} J.~Heinonen, {\it Lectures on Analysis on Metric
Spaces}, Springer-Verlag, 2001.

\bibitem {h-s} E.~Hewitt and K.~Stromberg, {\it Real and Abstract
Analysis}, Springer-Verlag, 1975.

\bibitem {h-w} W.~Hurewicz and H.~Wallman, {\it Dimension Theory},
Princeton University Press, 1941.

\bibitem {j1} B.~Jones, {\it A Class of Singular Integrals}, American
Journal of Mathematics {\bf 86} (1964), 441--462.

\bibitem {j2} F.~Jones, {\it Lebesgue Integration on Euclidean
Spaces}, Jones and Bartlett, 1993.

\bibitem {jl} J.-L.~Journ\'e, {\it Calder\'on--Zygmund Operators,
Pseudodifferential Operators, and the Cauchy Integral of Calder\'on},
Lecture Notes in Mathematics {\bf 994}, Springer-Verlag, 1983.

\bibitem {k} Y.~Katznelson, {\it An Introduction to Harmonic
Analysis}, 3rd edition, Cambridge University Press, 2004.

\bibitem {jk} J.~Kigami, {\it Analysis on Fractals}, Cambridge
University Press, 2001.

\bibitem {ak1} A.~Knapp, {\it Basic Real Analysis}, Birkh\"auser,
2005.

\bibitem {ak2} A.~Knapp, {\it Advanced Real Analysis}, Birkh\"auser,
2005.

\bibitem {sk1} S.~Krantz, {\it A Panorama of Harmonic Analysis},
Mathematical Association of America, 1999.

\bibitem {sk2} S.~Krantz, {\it Function Theory of Several Complex
Variables}, 2nd edition, AMS Chelsea Publishing, 2001.

\bibitem {sk3} S.~Krantz, {\it Real Analysis and Foundations}, 2nd
edition, Chapman \& Hall / CRC, 2005.

\bibitem {k-p} S.~Krantz and H.~Parks, {\it The Geometry of Domains in
Space}, Birkh\"auser, 1999.

\bibitem {sl} S.~Lang, {\it Real and Functional Analysis}, 3rd
edition, Springer-Verlag, 1993.

\bibitem {lax} P.~Lax, {\it Functional Analysis}, Wiley, 2002.

\bibitem {m-s-1} R.~Mac\'{\i}as and C.~Segovia, {\it Lipschitz
functions on spaces of homogeneous type}, Advances in Mathematics {\bf
33} (1979), 257--270.

\bibitem {m-s-2} R.~Mac\'{\i}as and C.~Segovia, {\it A decomposition
into atoms of distributions on spaces of homogeneous type}, Advances
in Mathematics {\bf 33} (1979), 271--309.

\bibitem {rm} R.~Ma\~n\'e, {\it Ergodic Theory and Differentiable
Dynamics}, translated from the Portuguese by S.~Levy, Springer-Verlag,
1987.

\bibitem {pm} P.~Mattila, {\it Geometry of Sets and Measures in
Euclidean Spaces: Fractals and Rectifiability}, Cambridge University
Press, 1995.

\bibitem {nad} M.~Nadkarni, {\it Basic Ergodic Theory}, 2nd edition,
Birkh\"auser, 1998.

\bibitem {pn} P.~Nicholls, {\it The Ergodic Theory of Discrete Groups},
Cambridge University Press, 1989.

\bibitem {kp} K.~Petersen, {\it Ergodic Theory}, Cambridge University
Press, 1989.

\bibitem {p} M.~Pollicott, {\it Lectures on Ergodic Theory and Pesin
Theory on Compact Manifolds}, Cambridge University Press, 1993.

\bibitem {p-y} M.~Pollicott and M.~Yuri, {\it Dynamical Systems and
Ergodic Theory}, Cambridge University Press, 1998.

\bibitem {pr} S.~Promislow, {\it A First Course in Functional
Analysis}, Wiley, 2008.

\bibitem {r} C.~Rickart, {\it General Theory of Banach Algebras}, van
Nostrand, 1960.

\bibitem {cr} C.~Rogers, {\it Hausdorff Measures}, with a foreword by
K.~Falconer, Cambridge University Press, 1998.

\bibitem {hr} H.~Royden, {\it Real Analysis}, 3rd edition, Macmillan,
1988.

\bibitem {ru1} W.~Rudin, {\it Principles of Mathematical Analysis},
3rd edition, McGraw-Hill, 1976.

\bibitem {ru2} W.~Rudin, {\it Function Theory in the Unit Ball of
${\bf C}^n$}, Springer-Verlag, 1980.

\bibitem {ru3} W.~Rudin, {\it Real and Complex Analysis}, 3rd edition,
McGraw-Hill, 1987.

\bibitem {ru4} W.~Rudin, {\it Functional Analysis}, 2nd edition,
McGraw-Hill, 1991.

\bibitem {dr} D.~Ruelle, {\it Dynamical Zeta Functions for Piecewise
Monotone Maps of the Interval}, American Mathematical Society, 1994.

\bibitem {r-y} B.~Rynne and M.~Youngson, {\it Linear Functional
Analysis}, 2nd edition, Springer-Verlag, 2008.

\bibitem {cs} C.~Sadosky, {\it Interpolation of Operators and Singular
Integrals: An Introduction to Harmonic Analysis}, Dekker, 1979.

\bibitem {ks} K.~Saxe, {\it Beginning Functional Analysis},
Springer-Verlag, 2002.

\bibitem {sch} M.~Schechter, {\it Principles of Functional Analysis},
2nd edition, American Mathematical Society, 2002.

\bibitem {slv} C.~Silva, {\it Invitation to Ergodic Theory}, American
Mathematical Society, 2008.

\bibitem {s1} Y.~Sinai, {\it Introduction to Ergodic Theory},
translated by V.~Scheffer, Princeton University Press, 1976.

\bibitem {s2} Y.~Sinai, {\it Topics in Ergodic Theory}, Princeton
University Press, 1994.

\bibitem {st1} E.~Stein, {\it Singular Integrals and Differentiability
Properties of Functions}, Princeton University Press, 1970.

\bibitem {st2} E.~Stein, {\it Topics in Harmonic Analysis Related to
the Littlewood--Paley Theory}, Princeton University Press, 1970.

\bibitem {st3} E.~Stein, {\it Boundary Behavior of Holomorphic
Functions of Several Complex Variables}, Princeton University Press,
1972.

\bibitem {st4} E.~Stein, {\it Harmonic Analysis: Real-Variable
Methods, Orthogonality, and Oscillatory Integrals}, with the
assistance of T.~Murphy, Princeton University Press, 1993.

\bibitem {st-sh-1} E.~Stein and R.~Shakarchi, {\it Fourier Analysis:
An Introduction}, Princeton University Press, 2003.

\bibitem {st-sh-2} E.~Stein and R.~Shakarchi, {\it Complex Analysis},
Princeton University Press, 2003.

\bibitem {st-sh-3} E.~Stein and R.~Shakarchi, {\it Real Analysis:
Measure Theory, Integration, and Hilbert Spaces}, Princeton University
Press, 2005.

\bibitem {s-w} E.~Stein and G.~Weiss, {\it Introduction to Fourier
Analysis on Euclidean Spaces}, Princeton University Press, 1971.

\bibitem {rs1} R.~Strichartz, {\it The Way of Analysis}, Jones and
Bartlett, 1995.

\bibitem {rs2} R.~Strichartz, {\it Differential Equations on
Fractals}, Princeton University Press, 2006.

\bibitem {str} K.~Stromberg, {\it Introduction to Classical Real
Analysis}, Wadsworth, 1981.

\bibitem {ds} D.~Stroock, {\it A Concise Introduction to the Theory of
Integration}, 3rd edition, Birkh\"auser, 1999.

\bibitem {t} M.~Taibleson, {\it Fourier Analysis on Local Fields},
Princeton University Press, 1975.

\bibitem {t1} A.~Torchinsky, {\it Real Variables}, Addison-Wesley,
1988.

\bibitem {t2} A.~Torchinsky, {\it Real-Variable Methods in Harmonic
Analysis}, Dover, 2004.

\bibitem {v-s-c} N.~Varopoulos, L.~Saloff-Coste, and T.~Coulhon, {\it
Analysis and Geometry on Groups}, Cambridge University Press, 1992.

\bibitem {w1} P.~Walters, {\it Ergodic Theory --- Introductory
Lectures}, Lecture Notes in Mathematics {\bf 458}, Springer-Verlag,
1975.

\bibitem {w2} P.~Walters, {\it An Introduction to Ergodic Theory},
Springer-Verlag, 1982.

\bibitem {w-z} R.~Wheeden and A.~Zygmund, {\it Measure and Integral:
An Introduction to Real Analysis}, Dekker, 1977.

\bibitem {rz} R.~Zimmer, {\it Ergodic Theory and Semisimple Groups},
Birkh\"auser, 1984.

\bibitem {z} A.~Zygmund, {\it Trigonometric Series}, Volumes I and II,
3rd edition, with a foreword by R.~Fefferman, Cambridge University
Press, 2002.





\end{thebibliography}
\end{document}